\documentclass{gtart}
\usepackage{amssymb, enumerate, epsfig}
\def\<{\langle}
\def\>{\rangle}

\def\C{{\mathbb C}}
\def\N{{\mathbb N}}
\def\Z{{\mathbb Z}}
\def\R{{\mathbb R}}

\def\pA{pseudo-Anosov}
\def\BnD{B_n / \langle \Delta^2 \rangle}
\def\tld{\widetilde}
\renewcommand{\phi}{\varphi}
\newtheorem{lemma}{Lemma}[section]
\newtheorem{theorem}[lemma]{Theorem}
\newtheorem{prop}[lemma]{Proposition}
\newtheorem{cor}[lemma]{Corollary}
\newtheorem{definition}[lemma]{Definition}

\newtheorem{remark}[lemma]{Remark}
\newtheorem{example}[lemma]{Example}

\begin{document}
\title{On the structure of the centralizer of a braid}
\author{Juan Gonz\'alez-Meneses\footnote{Partially supported by MCYT,
BFM2001-3207 and FEDER.} and Bert Wiest}
\address{Dpto.~de Matem\'atica Aplicada I, E.T.S.~Arquitectura, Universidad
de Sevilla, Avda.~Reina Mercedes, 2. 41012 Sevilla, Spain;
{\tt meneses@us.es}\\
and\\
IRMAR (UMR 6625 du CNRS), Universit\'e de Rennes 1, Campus de Beaulieu,
\\ 35042 Rennes cedex, France; {\tt bertold.wiest@math.univ-rennes1.fr}}

\begin{abstract}
The mixed braid groups are the subgroups of Artin braid groups whose elements preserve a
given partition of the base points. We prove that the centralizer of any braid can be
expressed in terms of semidirect and direct products of mixed braid groups. Then we
construct a generating set of the centralizer of any braid on $n$ strands, which has at most
$\frac{k(k+1)}{2}$ elements if $n=2k$, and at most $\frac{k(k+3)}{2}$ elements if $n=2k+1$.
These bounds are shown to be sharp, due to work of N.V.Ivanov and of S.J.Lee. Finally, we
describe how one can explicitly compute this generating set.
\end{abstract}
\makeshorttitle


\section{Introduction and statement of the results}\label{S:IN}

In 1971, Makanin \cite{Makanin} gave an algorithm for computing a
generating set of the centralizer $Z(\beta)$ of any given element $\beta$
of the $n$-string braid group $B_n$. His method, however, tends to yield
very large, and highly redundant generating sets.
One hint that much smaller generating sets could be found came from the experimental results
of Gonz\'alez-Meneses and Franco, which were obtained with a radically improved version of
Makanin's algorithm, based on new theoretical work \cite{GMFr}. Also, it has probably been
clear to specialists for a long time that Nielsen-Thurston theory could be used to improve
upon Makanin's results. However, there seems to be no such result in the literature, and the
aim of the present paper is to fill this gap.

Although our main interest was to compute, for any given $\beta\in B_n$,
a small generating set of $Z(\beta)$, we succeed in describing this
centralizer in terms of semidirect and direct products of
{\em mixed braid groups} (see~\cite{Manfredini,Orevkov}). These groups
are defined as follows: let $X=\{P_1,\ldots,P_n\}$ be the base points of
the braids in $B_n$. Given a partition ${\cal P}$ of $X$, the mixed braid
group $B_{\cal P}$ consists of those braids whose associated permutation
preserves each coset of ${\cal P}$.

The well known classification of mapping classes of a punctured surface
into periodic, reducible and pseudo-Anosov ones, yields an analogous
classification for braids. If $\beta$ is reducible, then one can decompose
it, in a certain sense, into a {\em tubular braid} $\widehat{\beta}$, and
some {\em interior braids} $\beta_{[1]},\ldots, \beta_{[t]}$, all of
them having less than $n$ stands. The main result of this paper is the
following:

\begin{theorem}\label{T:main1}
Let $\beta\in B_n$. One has:
\begin{enumerate}
\item If $\beta$ is pseudo-Anosov , then $Z(\beta)\simeq \mathbb{Z}^2$.

\item If $\beta$ is periodic, then $Z(\beta)$ is either $B_n$ or isomorphic
to a braid group on an annulus.

\item If $\beta$ is reducible, then there exists a split exact sequence:
$$
1 \longrightarrow Z(\beta_{[1]})\times \cdots \times Z(\beta_{[t]})
\longrightarrow Z(\beta) \longrightarrow Z_0(\widehat{\beta})\longrightarrow 1,
$$
where $Z_0(\widehat{\beta})$ is a subgroup of $Z(\widehat{\beta})$,
isomorphic either to $\Z^2$ or to a mixed braid group.
\end{enumerate}
\end{theorem}

Notice that $\mathbb{Z}\simeq B_2=B_{\{\{1,2\}\}}$, also
$B_n=B_{\{\{1,\ldots,n\}\}}$, and finally the braid group over an annulus
on $k$ strands is isomorphic to $B_{\{\{1,\ldots,k\},\{k+1\}\}}\subset
B_{k+1}$. Hence all these groups can be seen as mixed braid groups. Then,
by recurrence on the number of strands we deduce the following:

\begin{cor}
 For every $\beta\in B_n$, the centralizer $Z(\beta)$ can be expressed in
terms of semidirect and direct products of mixed braid groups.
\end{cor}

Using the above structure we shall construct, for any braid $\beta\in B_n$,
a generating set of $Z(\beta)$ having very few elements. More precisely,
we obtain:

\begin{theorem}\label{T:main}
If $\beta\in B_n$, then the centralizer $Z(\beta)$ can
be generated by at most $\frac{k(k+1)}{2}$ elements if $n=2k$,
and at most $\frac{k(k+3)}{2}$ elements if $n=2k+1$.
\end{theorem}

We will present an example, communicated to us by S. J. Lee, showing that
the above bound is sharp.
That is, we will define, for every positive
integer $n$, a braid in $B_n$ whose centralizer cannot be generated by
less than $\frac{k(k+1)}{2}$ elements if $n=2k$, or less than
$\frac{k(k+3)}{2}$ elements if $n=2k+1$.
(The first to observe that the number of generators of the centralizer
may grow quadratically with the number of strands was N.V.Ivanov
\cite{IvanovTalk}.)

However, the above bound refers to the worst case, and one could be interested
in the minimal number of generators of a particular braid. We shall give
a generating set which is in some sense the smallest ``natural'' generating
set for the centralizer of a braid. However, we shall also give an example
that illustrates the difficulty of finding the absolutely minimum possible
number of generators.

Let us mention that, for the special case of reducible braids conjugated
to a generator $\sigma_i$, its centralizer has already been described
in~\cite{FRZ}.

The plan of the paper is as follows: in section \ref{S:NiTh} we set up
notation and some standard machinery, and give the mentioned example by
S.\ J.\ Lee. In section~\ref{S:per} we study $Z(\beta)$ in the case where
$\beta$ is periodic, section~\ref{S:pA} deals with the {\pA} case, and
section~\ref{S:red} the reducible one, which is the most involved. In
section~\ref{S:upperbound} we define a generating set which is no larger
than the stated upper bound. In section~\ref{S:minimal} we describe a
generating set which is as small as possible while still reflecting the
geometric structure of the Nielsen-Thurston decomposition. We also give an
example to show that by algebraic trickery, even smaller sets can be obtained.
Finally in section~\ref{S:alg} we discuss how the generating set that
we defined can be found algorithmically.


\section{Prerequisites from Nielsen-Thurston theory}\label{S:NiTh}

We denote by $D$ the closed disk of radius $2$ centered at $0$ in the
complex plane. For any $n\in \N$, the disk $D$, together with any choice
of $n$ distinct points in its interior, is denoted $D_n$, and the
distinguished points are called the {\sl punctures}. We shall use
different choices for the exact position of the punctures at different
times - they may be lined up on the real axis, or regularly distributed on
a circle of radius $1$, or again one of them may be in the centre while
the remaining $n-1$ are distributed over the circle of radius $1$. In most
instances, the position of the punctures is irrelevant, and we shall leave
it unspecified.

We recall that the braid group $B_n$ is the group of isotopy classes of
homeomorphisms fixing (pointwise) the boundary and permuting the punctures
of $D_n$. Here the isotopies must fix pointwise the boundary and the punctures.
Alternatively, $B_n$ could be defined as the group of isotopy
classes of disjoint movements of the punctures, starting and ending with
the configuration of $D_n$. Yet another definition of $B_n$ is as the set
of isotopy classes of braids with $n$ strings in the cylinder
$D\times [0,1]$, where the start and end points of the strings are exactly
the puncture points in $D_n\times \{0\}$ and $D_n\times \{1\}$. We shall
use all three points of view.

We shall often work with a certain quotient of the group $B_n$, rather
than with $B_n$ itself.
We recall that the center of $B_n$ is isomorphic to the integers,
and generated by the full twist $\Delta^2$ (where $\Delta$ is Garside's
half twist). Geometrically, the group projection
$B_n \to B_n/\langle \Delta^2 \rangle$ is given by smashing the boundary
curve of $D_n$ to a puncture, so that $\BnD$ is naturally a subgroup
of the mapping class group of the sphere with $n+1$ punctures.
In order to keep notation manageable, we shall use the same letters for
elements of the braid group $B_n$ and for their image in the quotient $\BnD$.
This abuse of notation should not cause confusion.

We say that an element $\beta\in B_n$
is {\sl periodic} if the element of $\BnD$ represented by
$\beta$ is of finite order. Equivalently, $\beta$ is periodic if there
exists a $k\in \N$  such that in $B_n$ we have that $\beta^k$ is equal
to some power of $\Delta^2$.

We say an element $\beta$ of $B_n$ is {\sl reducible} if there exists a
nonempty multicurve $C$ in $D_n$ (i.e.~a system of disjoint simple closed
curves in $D_n$, none of them isotopic to the boundary or enclosing a
single puncture)
which is stabilized by $\beta$, i.e.~such that $\beta(C)$ is isotopic to $C$.
Note that $\beta$ may permute different components of the multicurve $C$.

The following definition is taken from \cite{BLM} (see also \cite{Ivanov}).
To every reducible braid $\beta\in B_n$ one can associate a canonical
invariant multicurve: its {\sl canonical  reduction system}, which by
definition is the collection of all isotopy classes $c$ of simple closed
curves which have the following two properties:
firstly, $c$ must be stabilized by some power of $\beta$, and secondly
any simple closed curve which has non-zero geometric intersection number
with $c$ must {\it not} be stabilized by any power of $\beta$. For instance,
let us
consider the punctured disk $D_6$, where the 6 punctures are arranged
uniformly on the circle of radius $1$ around $0$. Then the
rotation of the punctures around the circle by an angle of $\frac{2\pi}{3}$
is a periodic element of $B_6$ (of period $3$), it is
also reducible (e.g.~the three simple closed curves encircling punctures
1 and 2, 3 and 4, and 5 and 6 respectively form an invariant multicurve),
but its canonical  reduction system is empty. This example, however, is
somewhat untypical: if a {\it non-periodic} braid is reducible, then
its canonical reduction system is nonempty (see \cite{Ivanov}).

If $C$ is an invariant multicurve of a reducible braid $\beta$,
then we define the {\sl tubular braid} induced by $\beta$ and $C$ to be the
braid on fewer strings obtained from $\beta$ by removing from $D_n$ all the
disks bounded by outermost curves of $C$, and collapsing each outermost
curve of $C$ to a puncture point. It should be stressed that this braid
is only defined up to conjugacy.

An alternative way to look at the same defintion is the following:
let us consider again $\beta$ as an isotopy class of $n$ disjoint strings
in $D\times [0,1]$ with extremal points at the puncture points of
$D_n\times\{0\}$ and $D_n\times \{1\}$, such that each disk
$D\times \{t\}$ intersects each string exactly once.
Now our picture can be completed by embedded cylinders in
$D\times [0,1]$ which are disjoint from each other and from the strings
of the braid, each of which intersects each disk $D\times \{t\}$ in exactly
one circle, and whose boundary components are exactly the outermost curves
of $C$ in $D\times \{0\}$ and $D\times \{1\}$. We can interpret the solid
cylinders bounded by these cylinders as ``fat strings'', and the resulting
braid with some fat strings is exactly the tubular braid defined above.

The {\em interior} braids induced by $\beta$ and $C$ are the braids on
fewer strings induced by $\beta$ at the interior of the discs bounded by
the outermost curves of $C$. They can be thought of as the braids `inside'
the tubes of the tubular braid. Therefore, for every reducible braid
$\beta$, and every invariant multicurve $C$, we can decompose $\beta$ into
one tubular braid and some interior braids -- as many as the number of
outermost curves in $C$.

Finally, we have the notion of a {\sl \pA} element of $B_n$, for which we
refer to \cite{FLP} or \cite{Ivanov}. Roughly speaking, $\beta\in B_n$
is {\pA} if it is represented by a homeomorphism of $D_n$ which preserves
two transverse measured foliations on $D_n$ (called the ``stable'' and the
``unstable'' foliation), while scaling the measure of the unstable one by
some factor $\lambda$ which is greater than 1, and the measure of the stable
one by $\frac{1}{\lambda}$.

Thurston's theorem \cite{Thclass,FLP} states that every irreducible
element of $B_n$ is either periodic or pseudo-Anosov.

We end this section with the promised example, due to S. J. Lee, that should
be helpful for understanding the relationship between the Nielsen-Thurston
decomposition and the centralizer subgroup of a braid $\beta\in B_n$.
This example was also found independently by N. V. Ivanov and
H. Hamidi-Tehrani~\cite{IvanovRecent}.

\begin{example}\label{E:purebraid} \rm
 Suppose that $n=2m$, and denote by $\sigma_i$ the standard generator of
$B_n$, in which the $i$th and the $(i+1)$st punctures permute their
positions in a clockwise sense. We define $\beta=\sigma_1 \sigma_3^2
\sigma_5^3 \cdots \sigma_{2m-1}^{m}$.

The canonical reduction system of $\beta$ consists of $m$ circles, the
$i$th one enclosing the punctures $2i-1$ and $2i$. The corresponding
tubular braid is the trivial braid of $B_m$, and the interior braids are,
respectively, $\sigma_1$, $\sigma_1^2$, \ldots, $\sigma_1^m$ (notice that
all of them are non-conjugate, since conjugate braids have the same
exponent sum).

Let $D_{(1)}, \ldots, D_{(m)}$ be the disks bounded by the above circles.
As we shall see, any braid that commutes with $\beta$ has to send each
disk $D_{(i)}$ to itself (since the interior braids are non-conjugate).
A generating set of the centralizer subgroup of $\beta$ is given by
\begin{enumerate}
\item[(i)] for each $i\in \{1,\ldots,m\}$, the braid $\sigma_{2i-1}$, whose
 support is contained in $D_{(i)}$,
\item[(ii)] any generating set for the pure braid group on $m$ strings $P_m$--
all the generators here act as the identity on
$D_{(1)}\cup\ldots\cup D_{(m)}$, and can be seen as a pure tubular braid on $m$
strings (tubes), where the $i$th tube starts and ends at $D_{(i)}$.
\end{enumerate}
It can be easily shown that, in this case, $Z(\beta)\simeq\Z^m\times P_m$.
The essential observation now is the following: it can be deduced by the
presentation given in~\cite{Bir}, that the abelianization of $P_m$ is
isomorphic to $\Z^{m(m-1)/2}$  (see also~\cite{Arnold}). Hence, the
abelianization of $Z(\beta)$ is isomorphic to $\Z^m\times \Z^{m(m-1)/2}$.
Therefore, at least $m+\frac{m(m-1)}{2}=\frac{m(m+1)}{2}$ generators are
needed for the centralizer of the braid $\beta$.

 The case when $n=2m+1$ is analogous. The braid proposed by S. J. Lee is:
$\beta=\sigma_2 \sigma_4^2 \sigma_6^3 \cdots \sigma_{2m}^m$. This time the
first strand is not enclosed by any curve of the canonical reduction
system of $\beta$, and one has: $Z(\beta)\simeq \Z^m \times P_{m+1}$.
Hence, in this case the minimal possible number of generators is
$m+\frac{m(m+1)}{2}=\frac{m(m+3)}{2}$.
\end{example}

By proving theorem~\ref{T:main}, we will show that the above examples are
the worst one can find.


\section{The periodic case}\label{S:per}

We have to start by describing the periodic elements of $B_n$.
In order to state this classification result, which is classical, we
need to define two braids.

If $D_n$ is the disk with $n$ punctures arranged regularly on the circle
of radius $1$, then the braid which we shall call $\delta_{(n)}$
is represented by a clockwise
movement of all punctures on this circle by an angle $\frac{2\pi}{n}$.
If no confusion is possible, we shall simply write $\delta$, without
indicating the number of strands (note that this braid is the Garside
element of the Birman-Ko-Lee structure of $B_n$~\cite{BKL}).

Similarly, if we think of $D_{n}$ as having one puncture in the centre, and
$n-1$ punctures arranged circularly around it, then we define
$\gamma_{(n)}\in B_{n}$ to
be the braid given by a circular movement of the $n-1$ punctures by an
angle of $\frac{2\pi}{n-1}$, while leaving the central puncture fixed.
Again, for simplicity we shall often only write $\gamma$ instead of
$\gamma_{(n)}$.

The result that classifies periodic braids, which is due to
Eilenberg \cite{Eilenberg} and de~Ker\'ekj\'art\'o \cite{deK}
(see~\cite{coko} for a modern exposition) is:

\begin{lemma}\label{L:per_gamma_delta}
Every periodic braid in $B_n$ is conjugate to a power of $\delta_{(n)}$ or
$\gamma_{(n)}$.
\end{lemma}

\begin{figure}[htb]
\centerline{
\epsfbox{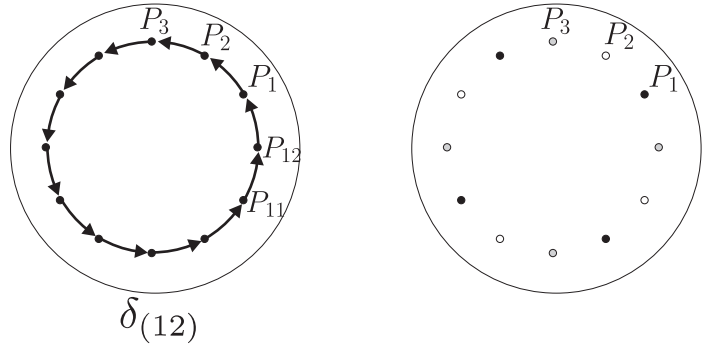}
}
\caption{The braid $\delta\in B_{12}$, and the three orbits of
$\delta^9$ (in black, white and grey).}
\label{F:delta_orbits}
\end{figure}

Thus we only need to consider the centralizer subgroups of $\delta_{(n)}^k$
and $\gamma_{(n)}^k$ for all $n, k\in \Z$, since the centralizers of
conjugate elements are isomorphic by an inner automorphism of $B_n$. This
problem has been solved by Bessis, Digne and Michel~\cite{BDM}, on the
wider context of complex reflexion groups. We shall explain their result in
the particular case of braid groups:

We suppose first that $\beta=\delta_{(n)}^k$ where, without loss of
generality, $k\geqslant 0$. Let $d=\gcd (n,k)$. For $u=1,\ldots,n$,
we will denote $P_u=e^{i 2\pi u/n}$ the punctures of $D_n$,
so $\beta=\delta_{(n)}^k$ sends $P_u$ to $P_{u+k}$ for every $u$
(the indices are taken modulo $n$). Hence the permutation induced
by $\beta$ has $d$ orbits (cycles) of length $r=\frac{n}{d}$, that we
denote by ${\cal C}_1,\ldots, {\cal C}_d$. See in figure~\ref{F:delta_orbits}
an example where $n=12$, $k=9$, $d=3$ and $r=4$: the braid $\delta_{(12)}$
and the three orbits of $\delta_{(12)}^9$.

\begin{figure}[htb]
\centerline{ \epsfbox{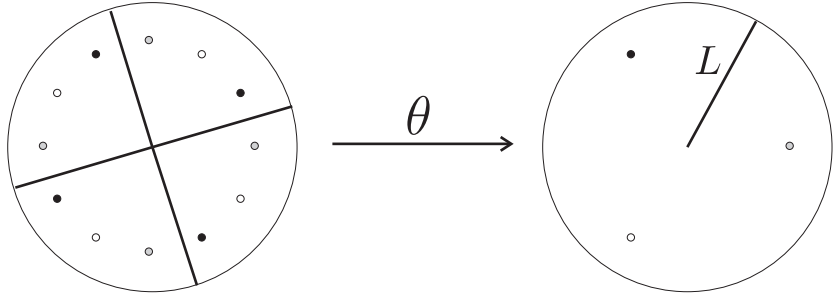} } \caption{The covering map $\theta=\theta_4$ associated to
$\delta_{(12)}^9$. } \label{F:theta}
\end{figure}

If $r>1$ (that is if $d<n$), consider the once punctured disc
$D^*=D\backslash \{0\}$, and the $r-$sheeted covering $\theta=\theta_r\co
D^* \rightarrow D^*$ defined by $\theta(ae^{it})=ae^{itr}=ae^{itn/d}$. The
orbits ${\cal C}_1,\ldots,{\cal C}_d$ are sent by $\theta$ to the points
$Q_1,\ldots,Q_d$, where $Q_u=e^{i2\pi u/d}$. If we consider the half-line
$L=\{ae^{i\pi/d},\; a\in ]0,2]\}$
(notice that $L$ passes between $Q_d$ and $Q_1$), then $D^*\backslash L$
is a fundamental region for $\theta$ (see figure~\ref{F:theta}).

Now notice that every braid in $B_d(D^*)$ can be lifted, by $\theta^{-1}$,
to a braid in $B_n$ in a natural way. The resulting braid is a $\frac{2\pi
d}{n}$-symmetric braid, that is, it is invariant under a rotation by an
angle of $\frac{2\pi d}{n}$. But then it is also invariant under a
rotation of angle $\frac{2\pi k}{n}$; in other words, the resulting braid
commutes with $\beta$. Hence we have a natural homomorphism: $\theta^*\co
B_d(D^*) \rightarrow B_n$ whose image is contained in $Z(\beta)$. Then one
has

\begin{theorem}[\cite{BDM}]\label{T:centrdelta}
The natural homomorphism
$\theta^*\co B_d(D^*)\rightarrow Z(\delta_{(n)}^k)$ is an isomorphism.
\end{theorem}

In other words, every element in the centralizer of $\beta=\delta_{(n)}^k$
can be seen (via $\theta$) as a braid on a once punctured disc, that is, a
braid on an annulus. Notice that if $r=1$ (that is, if $k$ is a multiple
of $n$), then $\beta$ is a power of $\delta_{(n)}^n=\Delta_{(n)}^2$. In this
case $\theta$ is the identity map, and the fundamental region is the whole
$D_n$. Hence the centralizer of $\beta$ is the whole $B_n$, as one should
expect.

\begin{figure}[htb]
\centerline{
\epsfbox{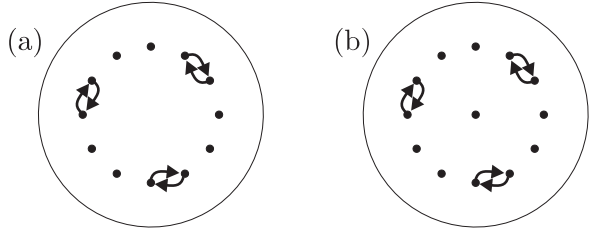} 
}
\caption{Generators $\theta_3^*(\sigma_1)$ and $\bar{\theta}_3^*(\sigma_1)$
of the centralisers of $\delta_{(12)}^4$ and $\gamma_{(13)}^4$.}
\label{F:gendelgam}
\end{figure}

Since we are interested in minimising the set of generators, we observe
that if $d=n$ (thus $r=1$), then $Z(\beta)=B_n$ is generated by two
elements, namely Artin's $\sigma_1$ and Birman-Ko-Lee's $\delta$. In a
similar way, if $1<d<n$, then the braid group $B_d(D^*)$ is generated by
just two elements, namely $\delta_{(n)}=\theta^*(\delta_{(d)})$ and the
braid $\theta^*(\sigma_1)$ shown in figure \ref{F:gendelgam}(a). Notice
that this case contains the above one, where $\theta^*$ is the identity.
Finally, if $d=1$ then $B_1(D^*)$ is cyclic. Thus we have:

\begin{prop}\label{P:centrdelta'}
If $k$ and $n$ are coprime, then $Z(\delta_{(n)}^k)$ is generated by a single
element, namely $\delta_{(n)}$. If, by contrast, $gcd(k,n)\geqslant 2$, then
$Z(\delta_{(n)}^k)$ is generated by two elements: $\delta_{(n)}$ and the braid
$\theta^*(\sigma_1)$.
\end{prop}

It is clear that the generating set given by proposition~\ref{P:centrdelta'}
is indeed minimal.
Next we study the centralizer of $\beta=\gamma_{(n)}^k$, still following the
work in \cite{BDM}.  This time we call $d=\gcd (n-1,k)$, and $r=(n-1)/d$.
If $d<n-1$, the above map $\theta$ induces a natural homomorphism
$\bar{\theta}^*=\bar{\theta}_r^*\co B_d(D^*)\rightarrow B_n$, where this
time the central point of $D$ is considered as a puncture. Hence, the central
strand of every braid coming from $B_d(D^*)$ is trivial. We observe that the
image of this homomorphism is contained in $Z(\beta)$, and in fact one has:

\begin{theorem}[\cite{BDM}]\label{T:centrdelta2}
The natural homomorphism
$\bar{\theta}_r^*\co B_d(D^*)\rightarrow Z(\gamma_{(n)}^k)$ is an isomorphism.
\end{theorem}

By contrast, if $d=n-1$, then $\beta$ is a power of $\gamma^{n-1}=\Delta^2$,
so $\overline{\theta}_r^*=1$, $Z(\beta)=B_n$ and everything works as above.
Hence we have

\begin{prop}\label{P:centrgamma}
If $k$ and $n-1$ are coprime, then $Z(\gamma_{(n)}^k)$ is generated by a
single element, namely $\gamma_{(n)}$. If, by contrast, $gcd(k,n-1)=d
\geqslant 2$, then $Z(\gamma_{(n)}^k)$ is generated by two elements:
$\gamma_{(n)}=\bar{\theta}^*(\delta_{(d)})$ and the braid
$\bar{\theta}^*(\sigma_1)$.
\end{prop}

See figure \ref{F:gendelgam}(b) for an illustration of the braid
$\bar{\theta}^*(\sigma_1)$. We summarize all the results in this
section as follows:

\begin{cor}
  The centralizer of any periodic braid in $B_n$ either equals $B_n$ or is
  isomorphic to $B_d(D^*)$, for some $d<n$. In particular, it can be
  generated by at most two elements.
\end{cor}

We end with a result that will be helpful later:

\begin{cor}\label{C:addstring}
If $k$ is not a multiple of $n$, then
$Z(\delta_{(n)}^k)\cong Z(\gamma_{(n+1)}^k)$.
\end{cor}

\begin{proof}
 Both groups are isomorphic to $B_d(D^*)$, where $d=\gcd(n,k)$.
An actual isomorphism can be defined as follows: take any element
$\alpha\in Z(\delta_{(n)}^k)$, isotope it to make it
$\frac{2\pi k}{n}$-symmetric, and then add a
trivial strand based at the central point of $D_n$.
\end{proof}


\section{The pseudo-Anosov case}\label{S:pA}

\begin{prop}\label{P:pA}
If $\beta\in B_n$ is pseudo-Anosov, then the centralizer of $B_n$ is free
abelian and generated by two elements: some {\pA} $\alpha$ which has the
same stable and unstable projective measured foliation as $\beta$ (possibly
$\beta$ itself), and one periodic braid $\rho$ (a root of $\Delta^2$,
possibly $\Delta^2$ itself).
\end{prop}

We stress that the generating set promised by proposition \ref{P:pA} is
obviously minimal.
For proving this result, it is more convenient to think about the quotient
group $\BnD$. Since $\langle \Delta^2 \rangle$ is the center of $B_n$,
it is contained in the centralizer of any element. Hence the centralizer of
an element in $B_n$ is just the preimage of the centralizer of its
corresponding mapping class in $\BnD$. Thus, for the rest of this section,
we shall work in this quotient $\BnD$; we shall prove the following
result, from which proposition \ref{P:pA} will then be deduced:

\begin{prop}\label{P:pA'} If $\beta\in \BnD$ is \pA, then the centralizer
of $\beta$ is abelian, and is generated by some {\pA} $\alpha$ which has the
same stable and unstable projective measured foliation as $\beta$, and
possibly one element $\rho$ of finite order.
\end{prop}

\begin{proof}[Proof of proposition \ref{P:pA'}]
We start by observing that the {\pA} element $\beta$ cannot commute with
any reducible element $a\in \BnD$, except possibly with periodic ones --
thus all elements of $Z(\beta)\subset \BnD$ are either {\pA} or periodic.
To see this, let us assume that the canonical reduction system $C$ of $a$
is non-empty. Then the canonical reduction system of
$\beta^{-1} a\beta$ is $\beta(C)$. If it were true that $\beta^{-1} a\beta=a$,
then we would have $\beta(C)=C$, which is impossible since it is well known
that {\pA} homeomorphisms do not stabilise any curves or multicurves.
(This result is also a special case of corollary 7.13 of \cite{Ivanov}.)

Our next claim is that all pseudo-Anosov elements in $Z(\beta)$ have the
same stable and unstable projective measured foliations.
In order to prove this, we can apply Corollaries 7.15 and 8.4 of
\cite{Ivanov}: since the centralizer subgroup of $\beta$
is infinite and irreducible,
it follows that $Z(\beta)$ contains an infinte cyclic group as a subgroup
of finite index. It follows that if $a$ is any {\pA} element in the
centralizer of $\beta$, then there exist $k, k'\in \N$ such that
$a^k=\beta^{k'}$. Since all powers of a {\pA} element have the same stable
and unstable projective measured foliation, it follows that $a$ has
the same stable and unstable projective measured foliations as $\beta$, and
so do all {\pA} elements of $Z(\beta)\subseteq \BnD$.

Next we make an essential observation which only works for braid groups,
and does not generalize to mapping class groups of surfaces with no boundary,
or with more than two boundary components: all elements of
$\BnD$, regarded as a subgroup of the mapping class group of the $n+1$ times
punctured sphere, fix the puncture which came from collapsing the boundary
of $D_n$. Moreover, there are singular leaves of the stable and unstable
foliation of $\beta$ emanating from this puncture, at least one of each
(like for every other puncture). In the cyclic ordering around the puncture,
singular leaves of the stable and unstable foliation alternate. If an
element $a$ of $\BnD$ commutes with $\beta$, then the action of $a$ has to
preserve the projective stable
and unstable foliations. Thus in the cyclic ordering around our preferred
puncture, the action of $a$ can only induce a cyclic (possibly trivial)
permutation of the singular leaves (sending stable to stable, and unstable
to unstable leaves, nevertheless).

Now we see that an element $a$ of $Z(\beta)\subseteq \BnD$ is uniquely
determined by just two data: firstly the stretch factor $\lambda$ by which
its action on the unstable measured foliation of $\beta$ multiplies the
measure of that foliation. (This factor $\lambda$ equals $1$ if $a$ is
periodic, and belongs to the set $\R_+\setminus \{1\}$ if $a$ is \pA).
And secondly by the cyclic permutation of the leaves of the stable projective
foliation emanating from the distinguished puncture of the $n+1$ times
punctured sphere. Indeed, if $a$ and $b$ share both data, then $ab^{-1}$ has
stretch factor 1 (so it is periodic), and preserves the singular leaves.
Hence it is the identity in $\BnD$, so $a=b$.

This implies that the set of periodic elements of
$Z(\beta)$ forms a subgroup of $Z(\beta)$ which is either trivial or
isomorphic to $\Z / k\Z$, where $k$ is a divisor of the number of singular
leaves of the stable foliation emanating from the preferred puncture.
Any generator of this subgroup can play the r\^ole of our desired generator
$\rho$ of $Z(\beta)\subseteq \BnD$.

Now $\rho$ commutes with any other element in $Z(\beta)$, because their
commutator has stretch factor 1 and induces the trivial permutation of the
prongs around the preferred singularity.

Now notice that the stretch factor yields a multiplicative map
from $Z(\beta)$ to $\mathbb{R}^+$. But it is known that the set of
possible stretch factors for a given foliation is discrete (see~\cite{Ivanov}),
so the image of $Z(\beta)$ under this map must be a cyclic subgroup of $\R^+$.
Take an element $\alpha$ whose stretch factor $\lambda$ generates this group.
Then $\alpha$ is {\pA} and the stretch factor of any element in
$Z(\beta)$ must be a power of $\lambda$.

We now have that $\alpha$ and $\rho$ generate $Z(\beta)\in \BnD$, because any
element in $Z(\beta)$ can be multiplied by some power of $\alpha$ so as
to obtain an element with stretch factor $1$, i.e.~a power of $\rho$.

It follows that $Z(\beta)\subset \BnD$ is isomorphic to $\Z \times \Z/ k\Z$,
with generators $\alpha$ and $\rho$. This completes the proof of
proposition \ref{P:pA'}.
\end{proof}

\begin{proof}[Proof of proposition~\ref{P:pA}]
 By proposition~\ref{P:pA'}, $Z(\beta)\subset \BnD$ is isomorphic to $\Z
\times \Z/ k\Z$, with generators $\alpha$ and $\rho$. But then
$Z(\beta)\subset B_n$ is just the preimage of $Z(\beta)\subset \BnD$ under
the natural projection. Consider the subgroup
$\langle \rho \rangle \subset Z(\beta) \subset \BnD$. Its preimage is an
infinte cyclic group in $B_n$ that contains $\langle \Delta^2 \rangle$.
We can suppose (up to choosing an appropriate $\rho$), that the generator
of this cyclic group projects to $\rho$, so we call it $\rho$ as well.
Notice that $\rho$ is a root of $\Delta^2$, since $\Delta^2$ belongs to
$\langle \rho \rangle$. Then we choose an element in $B_n$ that projects to
$\alpha$, and we also call it $\alpha$. We must prove that in $B_n$ we still
have $Z(\beta) = \langle \alpha \rangle \times \langle \rho \rangle$.

But every element in $Z(\beta)\subset B_n$ can be written as $\alpha^k
\rho^l \Delta^{2m}$. Since $\Delta^2$ is a power of $\rho$, then
$\{\alpha,\rho\}$ is a set of generators of $Z(\beta)$. On the other hand,
the commutator of $\alpha$ and $\rho$ projects to the trivial mapping
class, hence it equals $\Delta^{2k}$ for some $k$. But the algebraic
number of crossings of the braid $\Delta^{2k}$ is $kn(n-1)$, while for the
commutator of any two elements this number is zero. Hence $k=0$, so
$\alpha$ and $\rho$ commute. Finally, it is well known that $B_n$ is
torsion-free, so $Z(\beta)$ is isomorphic to $\Z \times \Z$, as we wanted
to show.
\end{proof}


\section{The reducible case}\label{S:red}

It remains to study the centralizer of a non-periodic reducible braid
$\beta$. Recall that for every braid $\gamma$ one has
$Z(\gamma^{-1}\beta\gamma)= \gamma^{-1} Z(\beta) \gamma$.
Hence, in general we will not study $Z(\beta)$, but the centralizer of a
suitable conjugate of $\beta$, which will be easier to describe.
Throughout this section we shall think of the punctures of the disk $D_n$
as being lined up on the real axis.


\subsection{Reducible braids in regular form}\label{S:reg_form}

As we saw in section~\ref{S:NiTh}, if $\beta$ is a non-periodic reducible
element, then its canonical reduction system is nonempty. We denote by
$R'(\beta)$ the set of outermost curves in the canonical reduction system
of $\beta$. It is determined by $\beta$ up to isotopy fixing the
punctures. Since we can study any conjugate of $\beta$, we can suppose
that $R'(\beta)$ is a family of disjoint circles centered at the real axis,
with disjoint interiors, each one enclosing more than one and less than
$n$ punctures.

Notice that there could be punctures in $D_n$ not enclosed by any circle
in $R'(\beta)$. In order to simplify the notations below, we define the
system of curves $R(\beta)$ to contain exactly the curves of $R'(\beta)$,
plus one circle around each such puncture of $D_n$. These new circles are
called the degenerate circles of $R(\beta)$. We now have that every puncture
in $D_n$ is enclosed by exactly one circle in $R(\beta)$.

 Notice that $\beta$ preserves $R(\beta)$, but it could permute the circles.
We will suppose that this permutation has $t$ orbits (or cycles)
${\cal C}_1, \ldots, {\cal C}_t$. That is, ${\cal C}_i$ is a family
of circles $\{C_{i,1},\ldots,C_{i,r_i}\}\subset R(\beta)$ such that $\beta$
sends $C_{i,k}$ to $C_{i,k+1}$ (here the second index is taken modulo
$r_i$). Then one has $R(\beta)={\cal C}_1\cup \cdots \cup {\cal C}_t =$
$\{C_{1,1},\ldots,C_{1,r_1}\}\cup \cdots \cup \{C_{t,1},\ldots,C_{t,r_t}\}$.
If $m_i$ is the number of punctures inside $C_{i,k}$, for any $k$, then
$1\leqslant m_i <n$ and $m_1r_1+\cdots +m_tr_t=n$.

Let $\widehat{\beta}$ be the tubular braid induced by $\beta$ and
$R(\beta)$. Then $\widehat{\beta}\in B_m$, where
$m=r_1+\cdots+r_t$. For $i=1,\ldots, t$ and $k=1,\ldots, r_i$, let
$\beta_{i,k}$ be the braid induced by $\beta$ in the interior of
$C_{i,k}$. In other words, $\beta_{i,k}$ is the braid inside the
tube of $\widehat{\beta}$ which starts at $C_{i,k}$ and ends at
$C_{i,k+1}$. We will call the braids $\beta_{i,k}$ the {\em interior braids}
of $\beta$. Notice that the interior braids of each degenerate circle is
just a trivial braid on one string.

\begin{figure}[htb]
\centerline{
 \epsfbox{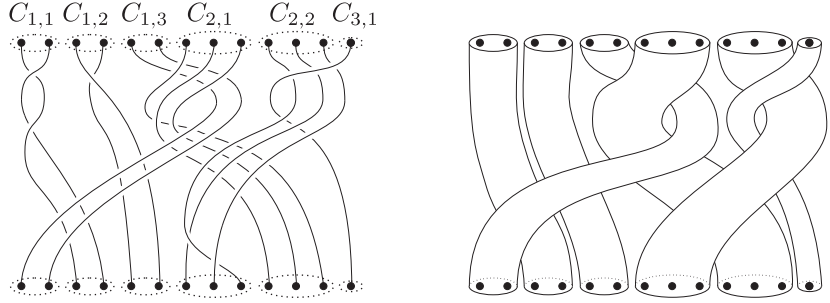} %
}
\caption{Example of a reducible braid $\beta$, and its corresponding tubular
braid $\widehat{\beta}$.}\label{F:tubular}
\end{figure}

In figure~\ref{F:tubular} we can see an example of a reducible braid
$\beta\in B_{13}$, and its corresponding tubular braid $\widehat{\beta}\in
B_6$. In this example we have three orbits, and the following data:
$r_1=3$, $r_2=2$, $r_3=1$; $m_1=2$, $m_2=3$, $m_3=1$,
$\beta_{1,1}=\sigma_1^2$, $\beta_{1,2}=\sigma_1^{-1}$, $\beta_{1,3}=1$,
$\beta_{2,1}=\sigma_1 \sigma_2$, $\beta_{2,2}=\sigma_1^{-1}\sigma_2$,
$\beta_{3,1}=1$ and $\widehat{\beta}=\sigma_3^2\sigma_2\sigma_1
\sigma_5^2\sigma_4$.

It would be desirable for $\beta$ to have its interior braids as simple
as possible, in order to study its centralizer. We propose the following:

\begin{definition}
Let $\beta\in B_n$ be a non-periodic reducible braid. Then $\beta$ will be
said to be in {\em regular form} if (using the notation introduced above) it
satisfies the following conditions:
\begin{enumerate}


 \item The only non-trivial interior braids in $\beta$ are
$\beta_{1,r_1}, \beta_{2,r_2},\ldots, \beta_{t,r_t}$ -- we shall denote these
braids by $\beta_{[1]}, \beta_{[2]}, \ldots, \beta_{[t]}$.

 \item For $i,j\in \{1,\ldots,t\}$, if $\beta_{[i]}$ and $\beta_{[j]}$ are
conjugate, then $\beta_{[i]}=\beta_{[j]}$.

\end{enumerate}
\end{definition}

 Hence, if $\beta$ is in regular form, there is at most one non-trivial
interior braid for each orbit, and any two interior braids are either equal or
non-conjugate. Fortunately, one can conjugate every non-periodic reducible
braid $\beta$ to another one in regular form, as we are going to see.

First, consider the subgroup $B_{R(\beta)}\subset B_n$ consisting of those
braids preserving $R(\beta)$. For $\alpha\in B_{R(\beta)}$, we can consider
the tubular braid $\widehat{\alpha}$ induced by $\alpha$ and $R(\beta)$.
Every $\alpha\in B_{R(\beta)}$ is completely determined by $\widehat{\alpha}$
and its interior braids $\alpha_{i,k}$, for
$i=1,\ldots t$ and $k=1,\ldots, r_i$.

 Now consider, in $\beta$, an orbit ${\cal C}_i =
\{C_{i,1},\ldots,C_{i,r_i}\}$ and the interior braids
$\beta_{i,1},\ldots, \beta_{i,r_i}\in B_{m_i}$. We define $\alpha\in
B_{R(\beta)}$ as follows: $\widehat{\alpha}$ is trivial, $\alpha_{j,k}=1$
if $j\neq i$, and $\alpha_{i,k}=\beta_{i,k}\beta_{i,k+1}\cdots
\beta_{i,r_i}$. If we conjugate $\beta$ by $\alpha$, we obtain
$\beta'=\alpha^{-1}\beta \alpha$, which has the following properties:
\begin{itemize} \item $\widehat{\beta}'=\widehat{\beta}$.

\item $\beta'_{j,k}=\beta_{j,k}$, for $j\neq i$.

\item $\beta'_{i,k}= (\alpha_{i,k})^{-1} \beta_{i,k} \alpha_{i,k+1} =
(\beta_{i,r_i}^{-1}\cdots \beta_{i,k}^{-1})(\beta_{i,k} \cdots
\beta_{i,r_i})=1$, for $k\neq r_i$.

\item $\beta'_{i,r_i}=(\alpha_{i,r_i})^{-1} \beta_{i,r_i} \alpha_{i,1} =
\beta_{i,r_i}^{-1} \beta_{i,r_i}(\beta_{i,1} \cdots \beta_{i,r_i}) =
\beta_{i,1} \cdots \beta_{i,r_i}$.

\end{itemize}

\begin{figure}[htb]
\centerline{
 \epsfbox{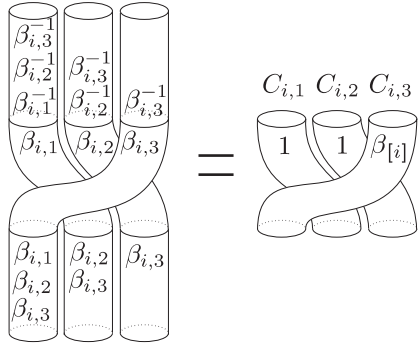} %
}
\caption{How to conjugate $\beta$ to simplify interior braids.}
\label{F:transfuse}
\end{figure}

In other words, if we conjugate $\beta$ by $\alpha$ we `transfuse' all the
interior braids in ${\cal C}_i$ to the last tube $C_{i,r_i}$, so
$\beta'_{i,r_i}$ becomes the only nontrivial interior braid in ${\cal
C}_i$. In figure~\ref{F:transfuse} we can see an example of such a
conjugation, where $\beta_{[i]}$ denotes the product $\beta_{i,1}\cdots
\beta_{i,r_i}$. We can now do the same for every $i=1,\ldots,t$.
Therefore, since we are interested in $\beta$ up to conjugacy, we can
suppose that $\beta_{i,k}=1$ if $k\neq r_i$ and denote
$\beta_{[i]}=\beta_{i,r_i}$, for every $i=1,\ldots,t$.

\begin{figure}[htb]
\centerline{
 \epsfbox{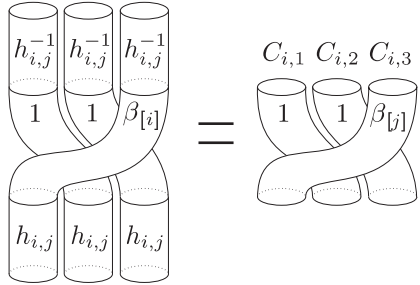} %
}
\caption{How to replace $\beta_{[i]}$ by $\beta_{[j]}$ if they are conjugate.}
\label{F:betai_betaj}
\end{figure}

Now suppose that some $\beta_{[i]}$ is conjugate to some $\beta_{[j]}$,
and let $h_{i,j}$ be a conjugating braid, that is, $h_{i,j}^{-1}
\beta_{[i]} h_{i,j}= \beta_{[j]}$. Consider the braid $\alpha\in
B_{R(\beta)}$ such that $\widehat{\alpha}=1$, $\;\alpha_{j,k}=1$ for
$j\neq i$ and $\alpha_{i,k}=h_{i,j}$ for every $k$. As we can see in
figure~\ref{F:betai_betaj}, if we conjugate $\beta$ by $\alpha$, then
$\beta_{[i]}$ is replaced by $\beta_{[j]}$. Therefore, we can assume that
for $i,j\in \{1,\ldots,t\}$, either $\beta_{[i]}=\beta_{[j]}$ or
$\beta_{[i]}$ and $\beta_{[j]}$ are not conjugate, and therefore we can
suppose that $\beta$ is in regular form.

\begin{figure}[htb]
\centerline{
 \epsfbox{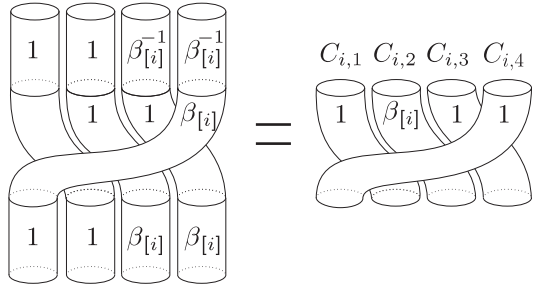} %
}
\caption{How to move $\beta_{[i]}$  from $C_{i,4}$ to $C_{i,2}$,
when $r_i=4$.}\label{F:mu_ik}
\end{figure}

 Notice that we have chosen to put $\beta_{[i]}$ into the tube starting at
$C_{i,r_i}$. But we can move it to any other tube of ${\cal C}_i$ if we
wish, by a suitable conjugation, and later on we will need to use this.
Hence we define, for $i\in \{1,\ldots,t\}$ and $k\in\{1,\ldots, r_i-1\}$,
a braid $\mu=\mu(i,k)$ that will `move' the interior braid $\beta_{[i]}$
to the tube $C_{i,k}$. This braid is defined as follows: the tubular braid
$\widehat{\mu}$ is trivial, and the interior braids are all trivial except
$\mu_{i,k+1}=\mu_{i,k+2}=\cdots =\mu_{i,r_i}=\beta_{[i]}$. We can see in
figure~\ref{F:mu_ik} how this works.

\subsection{Centralizer of a braid in regular form}

 We will now study the centralizer of $\beta$, assuming that $\beta$ is in
regular form. Recall that the only non-trivial interior braids of $\beta$
are denoted $\beta_{[1]},\ldots,\beta_{[t]}$, and that $\widehat{\beta}$
is the tubular braid associated to $\beta$ and $R(\beta)$. In this section
we will show that there is an exact sequence:
$$
  1 \rightarrow Z(\beta_{[1]})\times \cdots \times Z(\beta_{[t]})
\stackrel{g}{\longrightarrow}
Z(\beta) \stackrel{p}{\longrightarrow} Z_0(\widehat{\beta}) \rightarrow 1,
$$
where $Z_0(\widehat{\beta})$ is a subgroup of $Z(\widehat{\beta})$. Later
on we will see that this sequence splits.

For $i\in \{1,\ldots,t\}$, consider the centralizer $Z(\beta_{[i]})$ in
$B_{m_i}$. We define a map
$g_i\co Z(\beta_{[i]}) \rightarrow B_{R(\beta)}$ as follows:
given $\gamma\in Z(\beta_{[i]})$, $g_i(\gamma)$
is the braid $\alpha\in B_{R(\beta)}$ satisfying $\widehat{\alpha}=1$,
$\;\alpha_{j,k}=1$ for $j\neq i$, and $\alpha_{i,k}=\gamma$ for
$k=1,\ldots,r_i$. We need to show the following:

\begin{prop}\label{P:gi_injective}
 The map $g_i$ defined above is an injective homomorphism, and its image
is contained in $Z(\beta)$.
\end{prop}

\begin{proof}
%
The map $g_i$ is given by the diagonal homomorphism
$Z(\beta_{[i]})\to Z(\beta_{[i]}) \times \ldots \times Z(\beta_{[i]})$
($r_i$ factors), followed by the homomorphism induced by an inclusion of
$r_i$ copies of an $m_i$-times punctured disk into $r_i$ disjoint subdisks
(each containing $m_i$ punctures) of $D_n$. By the results of
\cite{ParisRolfsen} we can deduce that $g_i$ is indeed an injective
homomorphism.

It remains to show that for every $\gamma\in Z(\beta_{[i]})$ one has
$\alpha=g_i(\gamma)\in Z(\beta)$. Since $\widehat{\alpha}$ is trivial,
$\widehat{\alpha^{-1}\beta \alpha}= \widehat{\alpha}^{-1}\widehat{\beta}
\widehat{\alpha}=\widehat{\beta}$. So we just need to show that the
interior braids of $\alpha^{-1} \beta \alpha$ and $\beta$ coincide.
  For $j\neq i$, the braids $\alpha_{j,k}$ are trivial for every $k$, so
$\left(\alpha^{-1} \beta \alpha\right)_{j,k}=\beta_{j,k}$.
 Now, for $k\neq r_i$, one has $\left(\alpha^{-1} \beta \alpha\right)_{i,k} =
\alpha_{i,k}^{-1} \; \beta_{i,k}\;  \alpha_{i,k+1} =
\gamma^{-1} 1 \gamma = 1= \beta_{i,k}$.
 Finally, since $\gamma$ commutes with $\beta_{[i]}$, one has
$\left(\alpha^{-1} \beta \alpha\right)_{i,r_i} =
\alpha_{i,r_i}^{-1} \; \beta_{i,r_i}\;  \alpha_{i,1} =
\gamma^{-1} \beta_{[i]} \gamma = \beta_{[i]}= \beta_{i,r_i}$.
Therefore $\alpha^{-1}\beta \alpha = \beta$, so the image of $g_i$ is
contained in $Z(\beta)$.
\end{proof}

\begin{prop}\label{P:g_injective}
 The map  $\; g\co Z(\beta_{[1]})\times \cdots \times Z(\beta_{[t]})
\longrightarrow Z(\beta)$ defined by
$g(\gamma_1,\ldots,\gamma_t)=g_1(\gamma_1)\cdots g_t(\gamma_t)$
is an injective homomorphism.
\end{prop}

\begin{proof} Given $\gamma\in Z(\beta_{[i]})$, the only nontrivial
strands in $g_i(\gamma)$ are those inside the tubes
$C_{i,1},\ldots, C_{i,r_i}$.
Hence if $i\neq j$, $\gamma\in Z(\beta_{[i]})$ and
$\delta\in Z(\beta_{[j]})$, then $g_i(\gamma)$ and $g_j(\delta)$ commute.
Since every $g_i$ is a homomorphism, this shows that $g$ is also a
homomorphism. But we know by the previous proposition that $g_i$ is
injective for $i=1,\ldots,t$. Using an argument similar to the proof of
proposition \ref{P:gi_injective}, one can deduce that $g$ is also injective.
\end{proof}

Now we will relate $Z(\beta)$ and $Z(\widehat{\beta})$.
Every braid in $Z(\beta)$ preserves the canonical reduction system of
$\beta$ (see~\cite{Ivanov}), so it must preserve $R(\beta)$. That is,
$Z(\beta)\subset B_{R(\beta)}$. Let $p\co B_{R(\beta)}\rightarrow B_m$ be
the homomorphism which sends $\alpha$ to $\widehat{\alpha}$, the tubular
braid induced by $\alpha$ and $R(\beta)$. If we take $\alpha\in Z(\beta)$
then $\beta=\alpha^{-1}\beta \alpha$, so
$p(\beta)=p(\alpha^{-1}\beta \alpha)=p(\alpha)^{-1}p(\beta) p(\alpha)$.
Hence $p(\alpha)$ commutes with $p(\beta)=\widehat{\beta}$. Therefore, if
we restrict $p$ to $Z(\beta)$ we get
$\;p\co Z(\beta)\rightarrow Z(\widehat{\beta})$.

Unfortunately, neither $p\co B_{R(\beta)}\rightarrow B_m$ nor its restriction
$\;p\co Z(\beta)\rightarrow Z(\widehat{\beta})$ are surjective, but we shall
see that the
elements in the image of $p$ in either case can be easily characterised by
the permutation they induce. Notice that $p$ induces a bijection
$\widetilde{p}$ from $R(\beta)$ to $\{P_1,\ldots,P_m\}$, the punctures of
$D_m$. We denote by $\tau$ the inverse of $\widetilde{p}$.

\begin{definition}
 Let $\eta\in B_m$, and let $\pi_{\eta}$ be the permutation induced by
$\eta$ on the punctures of $D_m$. We say that $\pi_{\eta}$ is
{\em consistent with} $R(\beta)$ if, for $i=1,\ldots,m$,
$\;\tau(P_i)$ and $\tau(\pi_{\eta}(P_i))$ enclose the same number of punctures.
\end{definition}

\begin{prop}\label{P:cons_Rbeta}
An element $\eta\in B_m$ is in the image of
$\;p\co B_{R(\beta)}\rightarrow B_m$ if and only if
$\pi_{\eta}$ is consistent with $R(\beta)$.
\end{prop}

\begin{proof}
 If $\eta$ is in the image of $p$, let $\alpha\in B_{R(\beta)}$ with
$p(\alpha)=\eta$. Then, for every $i=1,\ldots,m$, $\;\tau(P_i)$ and
$\tau(\pi_{\eta}(P_i))$ are the top and bottom circles of a tube
determined by $\alpha$. Hence they must enclose the same number of
punctures (the number of strands inside the tube).

Conversely, suppose that $\pi_{\eta}$ is consistent with $R(\beta)$. Take
$i\in \{1,\ldots,m\}$ and suppose that $\tau(P_i)=C_{j,k}$. Then take the
$i$th strand of $\eta$ and consider it as a tube, enclosing the trivial
braid on $m_j$ strands. Do this for every $i=1,\ldots,m$. The resulting
braid, $\psi(\eta)$, is well defined since $\pi_{\eta}$ is consistent with
$R(\beta)$, and it belongs to $B_{R(\beta)}$. Moreover,
$p(\psi(\eta))=\eta$ by construction.
\end{proof}

The homomorphism $\psi$ introduced in this proof will play a prominent
r\^ole in what follows: if $\eta\in B_m$, then $\psi(\eta)$ is the braid
in $B_{R(\beta)}$ whose tubular braid equals $\eta$, and whose interior
braids are all trivial.

All the elements in $B_m$ that shall be considered from now on will have
permutations consistent with $R(\beta)$. Hence, by abuse of notation, we
will identify $C_{i,k}=\widetilde{p}(C_{i,k})$ and
${\cal C}_i=\widetilde{p}({\cal C}_i)$ if it does not lead to confusion.

We still need to characterise the elements in the image of
$\;p\co Z(\beta)\rightarrow Z(\widehat{\beta})$. We just know that their
permutations must be consistent with $R(\beta)$, but this is not
sufficient. Recall that the permutation induced by $\beta$ on the
components of $R(\beta)$ has $t$ orbits, ${\cal C}_1, \ldots {\cal C}_t$.
The key observation now is that every element $\alpha\in Z(\beta)$
preserves these orbits setwise, though it could permute them. Therefore,
for $i=1,\ldots,t$, one has $\alpha({\cal C}_i)={\cal C}_j$ for some $j$.
In the same way, for any $\eta\in Z(\widehat{\beta})$ one has
$\alpha({\cal C}_i)={\cal C}_j$ for some $j$.

\begin{lemma}\label{L:perm-conj}
Let $\alpha\in Z(\beta)$. If $\alpha({\cal C}_i)={\cal C}_j$ for some
$i,j\in \{1,\ldots, t\}$, then $\beta_{[i]}=\beta_{[j]}$.
\end{lemma}

\begin{proof}
Since $\alpha({\cal C}_i)={\cal C}_j$, the two orbits have the same length,
which we shall denote $r$; thus $r=r_i=r_j$.
Now $\beta^r$ is a braid that preserves $C_{i,k}$ and $C_{j,k}$ for every $k$,
and is such that $(\beta^r)_{i,k}=\beta_{[i]}$ and
$(\beta^r)_{j,k}=\beta_{[j]}$.
Now since $\alpha$ commutes with $\beta$, then it also commutes with
$\beta^{r}$. Suppose that $\alpha$ sends $C_{i,1}$ to $C_{j,k}$. Then
$\beta_{[j]}=(\beta^{r})_{j,k} = (\alpha^{-1}\beta^{r}\alpha)_{j,k} =
(\alpha_{i,1})^{-1} (\beta^{r})_{i,1} \alpha_{i,1} =
(\alpha_{i,1})^{-1} \beta_{[i]}\alpha_{i,1}.$
Therefore $\beta_{[i]}$ and $\beta_{[j]}$ are conjugate, and since $\beta$
is in regular form, $\beta_{[i]}=\beta_{[j]}$, as we wanted to prove.
\end{proof}

Lemma~\ref{L:perm-conj} imposes another condition for a braid in
$Z(\widehat{\beta})$ to be in $p(Z(\beta))$:

\begin{definition}
 Let $\eta\in Z(\widehat{\beta})$. We say that $\pi_{\eta}$ is
{\em consistent with} $\beta$ if it is consistent with $R(\beta)$ and,
furthermore, for every $i,j\in\{1,\ldots, t\}$ such that
$\eta({\cal C}_i)={\cal C}_j$, one has $\beta_{[i]}=\beta_{[j]}$.
\end{definition}

\begin{definition}
$Z_0(\widehat{\beta})$ is the subgroup of $Z(\widehat{\beta})$ consisting of
those elements whose permutation is consistent with $\beta$.
\end{definition}

Then lemma~\ref{L:perm-conj} can be restated as follows: If
$\alpha\in Z(\beta)$ then $p(\alpha)\in Z_0(\widehat{\beta})$. Moreover,
we can prove the following:

\begin{prop}\label{P:p_surjective}
 The homomorphism $p\co Z(\beta)\longrightarrow Z_0(\widehat{\beta})$ is
surjective.
\end{prop}

\begin{proof}
  Let $\eta\in Z_0(\widehat{\beta})$. We shall construct a preimage of $\eta$
under $p$ in two steps. Since $\pi_{\eta}$ is consistent with $\beta$
(thus with $R(\beta)$), we can, as a first step, consider the braid
$\psi(\eta)\in B_n$.
We then have $p(\psi(\eta))=\eta$; but $\psi(\eta)$ does not necessarily
commute with
$\beta$, since the interior braids of $\psi(\eta)^{-1}\beta \psi(\eta)$
could differ from those of $\beta$. Actually, since the interior braids of
$\psi(\eta)$ are all trivial, conjugating $\beta$ by $\psi(\eta)$ just
permutes the interior braids of $\beta$. More precisely, the braid
$\psi(\eta)^{-1} \beta \psi(\eta)$ equals $\beta$, except that, for each
$i\in\{1,\ldots,t\}$, it may not be the tube $C_{i,r_i}$ which contains the
nontrivial interior braid $\beta_{[i]}$, but some other tube from the
family $\mathcal{C}_i$. Our aim in the second step is thus to
fill the tubes of $\psi(\eta)$ with more suitable interior braids, in order
to obtain a braid that commutes with $\beta$.

 For every $i\in\{1,\ldots,t\}$, we know that $\psi(\eta)$ sends ${\cal
C}_i$ to some ${\cal C}_j$. Let $k_i\in \{1,\ldots,r_i\}$ be such that
$\psi(\eta)$ sends $C_{i,k_i}$ to $C_{j,r_j}$, and consider the braid
$\mu(i,k_i)$ defined at the end of Subsection~\ref{S:reg_form}. If we
conjugate $\beta$ by $\mu(i,k_i)$ we move $\beta_{[i]}$ from $C_{i,r_i}$
to $C_{i,k_i}$. If we further conjugate by $\psi(\eta)$, then
$\beta_{[i]}$ goes to $C_{j,r_j}$. But $\eta$ is consistent with $\beta$,
so $\beta_{[i]}=\beta_{[j]}$. Hence, the interior braids in ${\cal C}_j$
are preserved. We can do this for $i=1,\ldots,t$, so we obtain that the
braid
$$
\left(\prod_{i=1}^t{\mu(i,k_i)}\right)\psi(\eta)
$$
commutes with $\beta$ and its tubular braid is $\eta$, so it is in
$p^{-1}(\eta)\cap Z(\beta)$. This shows the result.
\end{proof}

We can finally bring together all the results in this section to state the
following:

\begin{theorem}\label{T:exactseq}
Let $\beta\in B_n$ be a non-periodic reducible braid in regular form.
Then the sequence
$$
1 \rightarrow Z(\beta_{[1]})\times \cdots \times Z(\beta_{[t]})
\stackrel{g}{\longrightarrow}
Z(\beta) \stackrel{p}{\longrightarrow}  Z_0(\widehat{\beta}) \rightarrow 1
$$
is exact.
\end{theorem}

\begin{proof}
  By proposition~\ref{P:g_injective} $g$ is injective, and by
proposition~\ref{P:p_surjective} $p$ is surjective. It just remains to
show that $\mbox{im}(g) =\ker(p)$.

By construction, every element in the image of $g$ induces a trivial
tubular braid, so $\mbox{im}(g) \subset \ker(p)$. Let then
$\alpha\in \ker(p)$, that is, $\widehat{\alpha}=1$. Since
$\alpha\in Z(\beta)$, we have $\alpha^{-1}\beta\alpha=\beta$, and since
$\beta_{i,k}=1$ for $k\neq r_i$, we must have
$\alpha_{i,k}^{-1} 1 \alpha_{i,k+1}=1$, so $\alpha_{i,k}=\alpha_{i,k+1}$ for
$k=1,\ldots,r_i-1$. Hence $\alpha_{i,1}=\alpha_{i,2}=\cdots = \alpha_{i,r_i}$
for every $i$. Moreover,  we have
$\beta_{[i]}=\beta_{i,r_i}=\alpha_{i,r_i}^{-1} \beta_{i,r_i} \alpha_{i,1}=
\alpha_{i,1}^{-1} \beta_{[i]} \alpha_{i,1}$, so
$\alpha_{i,1}\in Z(\beta_{[i]})$. Therefore,
$\alpha=g_1(\alpha_{1,1})g_2(\alpha_{2,1})\cdots g_t(\alpha_{t,1})=
g(\alpha_{1,1},\alpha_{2,1},\ldots,\alpha_{t,1})$. That is,
$\ker(p)\subset \mbox{im}(g)$.
\end{proof}


\subsection{Finding a section for $\mathbf{p}$}

 In this subsection we will prove that the exact sequence of
theorem~\ref{T:exactseq} splits. We recall that $\widehat{\beta}$ is
obtained from $\beta$ by collapsing the disks bounded by outermost curves
in the canonical reduction system of $\beta$ to single punctures. In
particular, the canonical reduction system of $\widehat{\beta}$ must be
empty. Hence, $\widehat{\beta}$ is either periodic or pseudo-Anosov. We will
distinguish these two cases, to define a multiplicative section for
$p$, but first we will show an easy particular case. Recall that a braid
is pure if it induces the trivial permutation of its base points.

\begin{prop}\label{P:section_trivial}
 If $\widehat{\beta}$ is pure, there is a homomorphism
$\:h\co Z_0(\widehat{\beta}) \rightarrow Z(\beta)$ such that $p\circ h=1$.
\end{prop}

\begin{proof}
  We shall prove that in this case, the homomorphism $\psi$ constructed
in the proof of proposition \ref{P:cons_Rbeta} is such a section.
  Let $\eta\in Z_0(\widehat{\beta})$. Since $\widehat{\beta}$ is pure,
${\cal C}_i=\{C_{i,1}\}$ for all $i$. Hence, if $\eta$ sends ${\cal C}_i$ to
${\cal C}_j$ then it sends the tube $C_{i,1}$ (containing $\beta_{[i]}$)
to the tube $C_{j,1}$ (containing $\beta_{[j]}=\beta_{[i]}$, since $\beta$
is in regular form). Therefore, filling every tube in $\eta$ with the
trivial braid, that is, defining $h(\eta)=\psi(\eta)$, yields indeed an
element of $Z(\beta)$.
\end{proof}

Next we study the general case, depending whether $\widehat{\beta}$ is
periodic or pseudo-Anosov.

\begin{prop}\label{P:section_periodic}
If $\widehat{\beta}$ is periodic, there is a homomorphism
$\:h\co Z_0(\widehat{\beta}) \rightarrow Z(\beta)$ such that $p\circ h=1$.
\end{prop}

\begin{proof}
Recall that we are studying $\beta$ up to conjugacy. This implies that we
can also study $\widehat{\beta}$ up to conjugacy since, for every
$\xi\in B_m$, if we conjugate $\beta$ by $\psi(\xi)$ we are conjugating
$\widehat{\beta}$ by $\xi$. Moreover, after conjugating by $\psi(\xi)$,
$\beta$ continues to be in regular form (up to renaming the circles in
$R(\beta)$). Therefore we can suppose, up to conjugacy, that
$\widehat{\beta}$ is a rigid rotation of the disc, that is, a power of
$\delta_{(m)}$ or $\gamma_{(m)}$.

Suppose first that $\widehat{\beta}=\delta^k_{(m)}$ for some $k$. We can
suppose that $k$ is not a multiple of $m$, since in that case
$\widehat\beta$ would be a power of $\Delta^2_{(m)}$, thus it would be pure,
and this case has already been studied in proposition~\ref{P:section_trivial}.
Recall the analysis of periodic braids in section~\ref{S:per}: the base
points $Q_1,\ldots, Q_m$ of $\widehat{\beta}$ will be evenly distributed
along a circle of radius 1 around 0. Let $d=\gcd(m,k)<m$ and $r=m/d$.
Then $\widehat\beta $ sends $Q_i$ to $Q_{i+k}$, and there are $d$ orbits
${\cal C}_1,\ldots,{\cal C}_d$ of length $r$. The orbit ${\cal C}_i$ will
contain the points $Q_u$ where $u\equiv i$ (mod $d$). Since we can choose
which tubes of $\beta $ contain the interior braids, we will suppose that
these are the tubes starting at $Q_{m-d+1}, Q_{m-d+2},\ldots, Q_{m}$, that
is, the last $d$ points of $D_m$.

We will consider now some line segments in $D$ which separate the points
$Q_1,\ldots,Q_m$ into $r$ sets of $d$ points. Let $L$ be the line segment
joining the origin with the border of $D$, passing between the points
$Q_{m-d}$ and $Q_{m-d+1}$, and let $L'$ be the segment passing between
$Q_m$ and $Q_1$. Notice that $L$ and $L'$ determine a sector which contains
the points $Q_{m-d+1},\ldots, Q_m$, corresponding to the tubes of $\beta$
with nontrivial interior braids. Let $\phi: \C \rightarrow \C $ be the
rotation around the origin by an angle of $2\pi k/m$ (the angle induced
by $\widehat\beta$), and denote $L_i=\phi^i(L)$. Since $\gcd(m,k)=d$, the
segments $L_0,\ldots,L_{r-1}$ divide $D$ into $m/d=r$ sectors, each
one of angle $2\pi/r$ and containing the points $Q_{id+1},\ldots, Q_{id+d}$
for some $i$. Take the smallest integer $e>0$ such that $\phi^e (L)=L'$.
Then one has $L_0=L$ and $L_e=L'$. We are interested in the union of segments
${\cal L}=L_1 \cup L_2\cup \cdots \cup L_e$ (see figure~\ref{F:L0_Le} for
an example).
\begin{figure}[htb]
\centerline{ \epsfbox{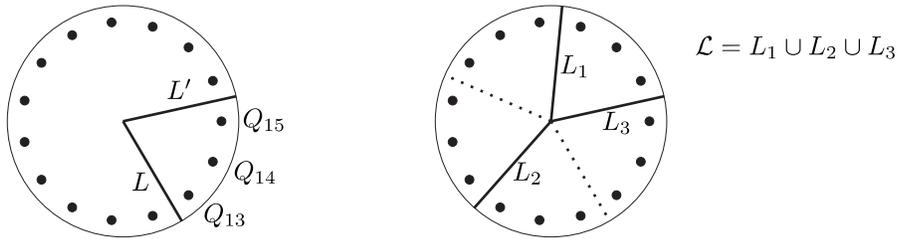} }
\caption{The segments $L$, $L'$, and the union of
segments ${\cal L}$, for $\widehat\beta=\delta^6\in B_{15}$.} \label{F:L0_Le}
\end{figure}

Let then  $\eta\in Z_0(\widehat{\beta})$.  In order to define $h(\eta)$, it
suffices to define its interior braids. This is done as follows: recall that,
since $\eta$ commutes with $\widehat{\beta}$, it can be isotoped to a
symmetric braid (with respect to the rotation $\phi $), so we take a
symmetric representative of $\eta$. For every base point $Q_i$ of
$\widehat{\beta}$ (corresponding to a circle $C_{j,u}$), consider the strand
of $\eta$ starting at $Q_i$ (the $i$th strand of $\eta$). Then we define the
interior braid $h(\eta)_{j,u}=(\beta_{[j]})^{L(\eta,i)}$, where
$L(\eta,i)\in \Z$ is the algebraic number of times that the $i$th strand of
$\eta$ crosses ${\cal L}$. This is well defined by theorem~\ref{T:centrdelta}
(if you take two distinct representatives of $\eta$ as a symmetric braid,
they are isotopic through symmetric braids, so the strands never touch the
origin and the intersection number $L(\eta,i)$ is preserved).

In other words, we define $h(\eta)$ as follows: we start with trivial
interior braids, and we follow the movement of the strands of $\eta$.
Each time a strand crosses a segment of ${\cal L}$ in the positive sense,
we multiply its interior braid by $\beta_{[j]}$ (where $j$ is the index
of the orbit ${\cal C}_j$ of that strand). And every time a strand crosses
${\cal L}$ in the negative sense, we multiply its interior braid by
$\beta_{[j]}^{-1}$.

We have thus defined a map
$h\co Z_0(\widehat{\beta})\rightarrow B_{R(\beta)}$. To show that
$h$ is a homomorphism, it suffices to see that the interior braids of
$\eta \xi$ are the product of those of $\eta$ and $\xi$, for
$\eta,\xi\in Z_0(\widehat{\beta})$. Suppose that the $i$th strand of
$\eta $ goes from $Q_i$ (corresponding to $C_{j,u}$) to $Q_{i'}$
(corresponding to $C_{j',u'}$). Hence $\eta $ sends ${\cal C}_j$ to
${\cal C}_{j'}$, and since $\eta\in Z_0(\widehat{\beta})$, it follows
that $\beta_{[j]}=\beta_{[j']}$. One also has, by definition,
$L(\eta\xi ,i)=L(\eta ,i)+ L(\xi ,i')$. Therefore $(\eta \xi)_{j,u} =
(\beta_{[j]})^{L(\eta\xi,i)} =
(\beta_{[j]})^{L(\eta,i)} (\beta_{[j]})^{L(\xi,i')} =
\eta_{j,u}\xi_{j',u'}$, so  $h$ is a homomorphism.

We must finally show that, with this definition, $h(\eta)\in Z(\beta )$,
for every $\eta \in Z_0(\widehat\beta )$. We will define first some special
braids. For every $i,j\in \{1,\ldots,d\}$ such that  $i<j$ and
$\beta_{[i]}=\beta_{[j]}$, define the symmetric braid $S_{i,j} =
S_{j,i}=\theta_r^*(\sigma_i \cdots \sigma_{j-2} \sigma_{j-1} \sigma _{j-2}
\cdots \sigma_i)$  (see figure~\ref{F:gendelgam} in section~\ref{S:per} to
recall the definition of $\theta_r^*$, and figure~\ref{F:S_13} here for an
example). The braid $S_{i,j}$ commutes with $\widehat\beta $ (since it is
symmetric), and it permutes the orbits ${\cal C}_i$ and ${\cal C}_j$,
preserving the others. Hence $S_{i,j}\in Z_0(\widehat\beta)$. Moreover, its
strands do not cross ${\cal L}$, so by definition of $h$ one has
$h(S_{i,j})=\psi(S_{i,j})$ (the interior braids are trivial).
\begin{figure}[htb]
\centerline{ \epsfbox{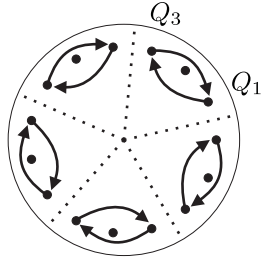} }
\caption{The braid $S_{1,3}$, for
$\widehat\beta=\delta^6\in B_{15}$ (assuming that $\beta_{[1]}=\beta_{[3]}$).}
\label{F:S_13}
\end{figure}

But $h(S_{i,j})$ commutes with $\beta$, since the only tubes it permutes
are those of the orbits ${\cal C}_i$ and ${\cal C}_j$; among these tubes,
the only two with non-trivial interior braids are exchanged, and their
corresponding interior braids are equal ($\beta_{[i]}=\beta_{[j]}$). Hence
the interior braids of $\beta$ are preserved by $\psi(S_{i,j})=h(S_{i,j})$,
so $h(S_{i,j})\in Z(\beta)$.

Take then an arbitrary $\eta\in Z_0(\widehat\beta )$. We must show that
$h(\eta) \in Z(\beta)$. Suppose that $\eta$ sends ${\cal C}_i$ to
${\cal C}_j$ for some $i,j$. Then $\beta_{[i]}=\beta_{[j]}$, so $S_{i,j}$
is defined, and the braid $\eta S_{i,j}$ preserves the orbit ${\cal C}_i$.
We can continue this way, until we obtain a braid
$\eta S_{i_1,j_1} \cdots S_{i_k,j_k}$ that commutes with $\widehat\beta $
and preserves every orbit ${\cal C}_i$, for $i=1,\ldots,d$.  Since
$h(S_{i,j})\in Z(\beta )$ for every $i,j$, and $h$ is a homomorphism, in
order to show that $h(\eta)\in Z(\beta)$ it suffices to show that $h(\eta
S_{i_1,j_1}\cdots S_{i_k,j_k})\in Z(\beta)$. Therefore, we can suppose that
$\eta$ preserves every orbit ${\cal C}_i$.

Denote $\alpha =h(\eta)$. We need to show that the interior braids of $\alpha^{-1} \beta
\alpha $ coincide with those of $\beta$. Since $\eta$ preserves all orbits, we will consider
just the tubes of ${\cal C}_1$, the other ones being analogous. Suppose that $\alpha$ sends
the circle $C_{1,u}$ to $C_{1,r}$. Then it must send $C_{1,v}$ to $C_{1,v-u}$ for every $v$
(the indices are taken modulo $r$).

We will identify the points $Q_1,\ldots,Q_m$ with their corresponding circles $C_{i,v}$. For
every $v=1,\ldots,r$, let $b_v$ be the strand of $\eta $ starting at $C_{1,v}$. Since $\eta$
is symmetric, we have $\phi (b_v)=b_{v+1}$. Suppose that $b_v$ crosses $t$ times the segment
$L_i$, where $i\in\{0,\ldots,r-1\}$. Then $b_{v+1}$ will cross $t$ times the segment
$\phi(L_i)=L_{i+1}$. Therefore, if $b_v$ crosses $l$ times ${\cal L}$, and if it crosses
$l_0$ times $L_0$ and $l_e$ times $L_e$, then $b_{v+1}$ crosses $l-l_e+l_0$ times ${\cal
L}$.

If $v\neq r$ and $v\neq u$, then $b_v$ neither starts nor ends at $C_{1,r}$. Then it crosses
$L_0$ and $L_e$ the same number of times. Hence, $b_v$ and $b_{v+1}$ cross ${\cal L}$ the
same number of times, say $l$. Therefore, if $v\neq r,u$, one has $(\alpha^{-1} \beta
\alpha)_{1,v-u}=(\alpha_{1,v})^{-1} \beta_{1,v} \alpha_{1,v+1} = \beta_{[1]}^{-l} \: 1 \:
\beta_{[1]}^l= 1 = \beta_{1,v-u}$.

If $u=v=r$, then $b_v$ starts at ends at $C_{1,r}$. Hence, as above, it crosses $L_0$ and
$L_e$ the same number of times, so $b_v=b_r$ and $b_{v+1}=b_1$ cross ${\cal L}$ the same
number of times, say $l$. We then have $(\alpha^{-1} \beta \alpha)_{1,v-u}=(\alpha^{-1}
\beta \alpha)_{1,r} = (\alpha_{1,r})^{-1} \beta_{1,r} \alpha_{1,1}= \beta_{[1]}^{-l}
\beta_{[1]} \beta_{[1]}^l = \beta_{[1]} = \beta_{1,r} = \beta_{1,v-u}$. Hence, if $u=r$, we
have already seen all the possible cases. We will then suppose that $u\neq r$.

If $v=r$, then $b_v$ starts (but does not end) at $C_{1,r}$. Hence, it crosses $L_e$ one
more time (in the positive sense) than it crosses $L_0$. Therefore, if $b_v=b_r$ crosses $l$
times ${\cal L}$, then $b_{v+1}=b_1$ crosses it $l-1$ times. One has: $(\alpha^{-1} \beta
\alpha)_{1,v-u}=(\alpha^{-1} \beta \alpha)_{1,r-u} = (\alpha_{1,r})^{-1} \beta_{1,r}
\alpha_{1,1}= \beta_{[1]}^{-l} \beta_{[1]} \beta_{[1]}^{l-1} = 1 =
\beta_{1,r-u}=\beta_{1,v-u}$.

Finally, if $v=u$ then $b_v$ ends (but does not start) at $C_{1,r}$. In this case, it
crosses $L_e$ one less time (in the positive sense) than it crosses $L_0$. Hence, if
$b_v=b_u$ crosses $l$ times ${\cal L}$, then $b_{v+1}=b_{u+1}$ crosses it $l+1$ times. One
then has: $(\alpha^{-1} \beta \alpha)_{1,v-u}=(\alpha^{-1} \beta \alpha)_{1,r} =
(\alpha_{1,u})^{-1} \beta_{1,u} \alpha_{1,u+1}= \beta_{[1]}^{-l} \: 1 \: \beta_{[1]}^{l+1} =
\beta_{[1]} = \beta_{1,r} = \beta_{1,v-u}$.

Therefore, in every possible case we have  $(\alpha^{-1} \beta \alpha)_{1,v-u}=
\beta_{1,v-u}$, for every $v$. This means that the interior braids of $(\alpha^{-1} \beta
\alpha)$ and of $\beta $ coincide, that is, $\alpha =h(\eta)$ commutes with $\beta $, as we
wanted to show.

This completes the proof of proposition \ref{P:section_periodic} in the case
$\widehat{\beta}=\delta^k_{(m)}$, and it only remains to deal with the case when
$\widehat\beta=\gamma_{(m)}^k$. As above, we can suppose that $k$ is not a multiple of
$m-1$, since in that case $\widehat\beta$ would be pure, and this case has already been
treated in proposition~\ref{P:section_trivial}. Hence, the only fixed point in the
permutation induced by $\widehat\beta $ is the origin. Therefore, every $\eta$ commuting
with $\widehat\beta $ must fix the origin. This means that, for every $\eta\in
Z_0(\widehat\beta )$, we can fill its central tube with the trivial braid, and the other
tubes in the same way as above (defining ${\cal L}$, and counting the number of times each
strand crosses ${\cal L}$). This defines a homomorphism $h: Z_0(\widehat\beta)\rightarrow
Z(\beta)$ which is a section of $p$. The proof is the same as above.
\end{proof}

It remains to study the case when $\widehat{\beta}$ is \pA.

\begin{prop}\label{P:section_pA}
If $\widehat{\beta}$ is pseudo-Anosov, then there is a homomorphism
$h\co Z_0(\widehat{\beta})\rightarrow Z(\beta)$ such that $p\circ h=1$.
\end{prop}

\begin{proof}
In this case, we know that $Z(\widehat{\beta})$ is a free abelian group of rank 2, generated
by a pseudo-Anosov and a periodic braid. Hence, $Z_0(\widehat{\beta})$ is an abelian group
of rank one or two. Notice that $\Delta_{(m)}^2\in Z_0(\widehat{\beta})$, because this braid
commutes with $\widehat{\beta}$ and because $\pi_{\Delta^2}$ is trivial, and thus consistent
with $\beta$. Hence $Z_0(\widehat{\beta})$ contains at least one periodic element. On the
other hand, $\widehat{\beta}$ belongs itself to $Z_0(\widehat{\beta})$, since
$\pi_{\widehat{\beta}}$ is clearly consistent with $\beta$. Hence in $Z_0(\widehat{\beta})$
there are also pseudo-Anosov braids. Since all powers of a periodic braid are periodic, and
all powers of a pseudo-Anosov braid are pseudo-Anosov, it follows that
$Z_0(\widehat{\beta})$ has in fact rank two.  More precisely, $Z_0(\widehat{\beta}) =
\langle \eta\rangle \times \langle \rho \rangle$, where $\eta$ is pseudo-Anosov and $\rho$
is periodic. In particular, we have
$\widehat\beta \in \langle \eta\rangle \times \langle \rho \rangle$, and the three braids
$\widehat{\beta}$, $\eta$ and $\rho$ are mutually commuting.

Our aim is to define two commuting braids $h(\rho)$ and $h(\eta)$ in
$Z(\beta)$ which are preimages of $\rho$ respectively $\eta$ under $p$.
The definition of $h(\rho)$ is very simple: we take an arbitrary preimage
of $\rho$ under $p$ -- this is possible since $p$ is surjective by
proposition \ref{P:p_surjective}. It remains to construct $h(\eta)$.

\begin{lemma}\label{L:pAlem1}
Suppose $\alpha\in B_{R(\beta)}$, that is, the braid $\alpha$ preserves the set of outermost
curves in the canonical reduction system of $\beta$. Suppose also that $\mu, \nu \in
Z(\widehat{\alpha})$. Suppose that $\iota_\mu\in B_{R(\beta)}$ is a braid with trivial tubes
(i.e.\ $\widehat{\iota}_\mu=1$) such that $\psi(\mu)\cdot\iota_\mu \in Z(\alpha)$. Finally,
suppose that $\mu$ and $\nu$ induce the same permutation. Then we have as well that
$\psi(\nu)\cdot\iota_\mu \in Z(\alpha)$.
\end{lemma}

In other words, if two tubular braids commute with $\widehat{\alpha}$, if they induce the
same permutation, and if some ``filling'' of one of them commutes even with $\alpha$, then
the same filling of the other will also commute with $\alpha$.

\begin{proof}[Proof of lemma \ref{L:pAlem1}] Conjugating $\alpha$ by
$\psi(\nu)\cdot\iota_\mu \in Z(\alpha)$ yields a certain braid $\alpha'$; we have to check
that $\alpha'=\alpha$. Firstly, we have an equality of tubular braids
$\widehat{\alpha'}=\widehat{\alpha}$, because $\nu$, the tubular braid of
$\psi(\nu)\cdot\iota_\mu$, commutes with $\widehat{\alpha}$. Moreover, since $\mu$ and $\nu$
induce the same permutations, we have for $i=1,\ldots, m$ that the $i$th tube of $\alpha'$
contains the same braid as the $i$th tube of
$(\psi(\mu)\cdot\iota_\mu)^{-1}\cdot\alpha\cdot(\psi(\mu)\cdot\iota_\mu)$. Since
$\psi(\mu)\cdot\iota_\mu$ commutes with $\alpha$, this is in turn the same as the $i$th tube
of $\alpha$. In summary, $\alpha$ and $\alpha'$ have the same tubular braids, and
corresponding tubes contain the same interior braids, which implies that $\alpha=\alpha'$.
\end{proof}

Next we have to think in detail about the orbit structure of
$\widehat{\beta}$. Let us choose arbitrarily a puncture $P$ of the disk
$D_m$ (on which $\widehat{\beta}$ acts), and let $O(\widehat{\beta},\rho)$ be
the orbit of that puncture under the action of the subgroup
$\langle \widehat{\beta}\rangle \times \langle \rho \rangle$ of
$Z_0(\widehat{\beta})$. Let $O(\widehat{\beta},\rho,\eta)$ be the orbit
of $P$ under the action of the group
$\langle\rho\rangle \times \langle \eta \rangle$ (note that this
group is also isomorphic to $\Z^2$, and contains $\widehat{\beta}$).

We are going to suppose without loss of generality that
$O(\widehat{\beta},\rho,\eta)$ contains {\it all} punctures of $D_m$,
and we shall specify how the tubes of $\eta$ corresponding to this orbit
shall be filled -- indeed, if there are other orbits, then these can be
treated in same way, independently.

{\bf Special case: } Let us start by considering the simpler special case that
$O(\widehat{\beta},\rho,\eta) = O(\widehat{\beta},\rho)$, i.e.\ that the
action of $\eta$ preserves the $(\widehat{\beta},\rho)$-orbit. In this case
we have
\begin{lemma}\label{L:pAlem2}
There exist integers $k$ and $l$ such that $\eta$ and
$\widehat{\beta}^k\cdot \rho^l$ induce the same permutations on
$O(\widehat{\beta},\rho)$.
\end{lemma}
\begin{proof}[Proof of lemma \ref{L:pAlem2}] One can choose $k$ and $l$ such
that $\widehat{\beta}^k \rho^l(P) = \eta(P)$, simply because $\eta(P)$ is in
the orbit of $P$ under the action of $\widehat{\beta}$ and $\rho$. Now if
$P'$ is another point in the orbit, then
$P'=\widehat{\beta}^\kappa \rho^\lambda(P)$ for some $\kappa, \lambda \in \Z$.
Since $\widehat{\beta}, \rho$, and $\eta$ are mutually commuting, we get
$\widehat{\beta}^k \rho^l(P') =
\widehat{\beta}^k \rho^l(\widehat{\beta}^\kappa \rho^\lambda(P)) =
\widehat{\beta}^\kappa\rho^\lambda(\widehat{\beta}^k \rho^l(P))=
\widehat{\beta}^\kappa\rho^\lambda(\eta(P))=
\eta(\widehat{\beta}^\kappa \rho^\lambda(P)) = \eta(P')$.
\end{proof}

We already know a nice preimage of $\widehat{\beta}^k\cdot \rho^l$ under $p$:
the braid $\beta^k\cdot h(\rho)^l$ belongs to $Z(\beta)$, because both
$\beta$ and $h(\rho)$ do. This braid can be reexpressed as
$\psi(\widehat{\beta}^k \rho^l)\cdot\iota$, where $\iota$ is
some braid in $B_{R(\beta)}$ with $\widehat{\iota}=1$. (That is, we define
$\iota$ to consist of the interior braids of the tubes of
$\beta^k\cdot h(\rho)^l$).

Now we define our filling of $\eta$ by $h(\eta):=\psi(\eta)\cdot\iota$. By
lemma \ref{L:pAlem1} we have that indeed $\psi(\eta)\cdot\iota\in Z(\beta)$.
In order to see that $\psi(\eta)\cdot\iota$ lies also in the centralizer of
$h(\rho)$ one can use a very similar argument. Explicitly, both $\eta$
and $\widehat{\beta}^k\rho^l$ lie in the centralizer of $\rho$, and they
induce the same permutation of the punctures. Moreover,
$\psi(\widehat{\beta}^k \rho^l) \cdot\iota = \beta^k h(\rho)^l\in Z(h(\rho))$.
By lemma \ref{L:pAlem1} we conclude again that
$\psi(\eta)\cdot\iota \in Z(h(\rho))$, also.

{\bf General case: } In the case where $\eta$ does not preserve
$O(\widehat{\beta},\rho)$, the strategy is to work not with $\beta$ itself
but with a certain conjugate of $\beta$. The details are as follows.
We have a finite number of disjoint $(\widehat{\beta},\rho)$-orbits in
$O(\widehat{\beta},\rho,\eta)$, and since $\eta$ commutes with
$\widehat{\beta}$ and $\rho$, the action of $\eta$ permutes these
orbits cyclically:
$$
O(\widehat{\beta},\rho)\stackrel{\eta-{\rm action}}{\longrightarrow}
\eta(O(\widehat{\beta},\rho)) \stackrel{\eta-{\rm action}}{\longrightarrow}
\ldots \stackrel{\eta-{\rm action}}{\longrightarrow}
\eta^s(O(\widehat{\beta},\rho))=O(\widehat{\beta},\rho).
$$
Let us denote $\widehat{\beta}_*, \rho_*$ and $\eta^s_*$ the braids
which are obtained from $\widehat{\beta}, \rho$ and $\eta^s$ by retaining
only the strands corresponding to $O(\widehat{\beta},\rho)$, and forgetting
the strands corresponding to all other $(\widehat{\beta},\rho)$-orbits.
Similarly, let $\beta_*$ be the corresponding restriction of $\beta$.
Our first aim is to fill the tubes of $\rho_*$ and $\eta^s_*$ so as to obtain
commuting braids in $Z(\beta_*)$. This can be done as in the ``special case'':
for $\rho_*$ we choose any filling in $Z(\beta_*)$, and for $\eta^s_*$ there
exists a braid $\iota_*$ with trivial tubes such that
$\psi(\eta_*^s) \cdot \iota_*$ commutes with $\beta_*$ and the filling of
$\rho_*$.

We have succeeded in finding a filling of certain tubes of $\eta^s$, but
not yet of $\eta$ itself. Also, we have so far only filled the tubes
of $\rho$ which correspond to $O(\widehat{\beta},\rho)$, but not yet
those in the $\eta$-translates of this orbit. We first notice that the
$\eta$-action sends $\widehat{\beta}$-orbits to $\widehat{\beta}$-orbits,
and that in each $\widehat{\beta}$-orbit there is exactly one tube whose
preimage in $\beta$ contains a nontrivial braid (the same for all
$\widehat{\beta}$-orbits), and all other tubes are filled with
a trivial braid. Thus, up to cyclically changing the numbering of the
orbits of each tube of $\widehat{\beta}$, we may assume that the
$\eta$-action sends each tube of $\widehat{\beta}$ in
$O(\widehat{\beta},\rho)$ to a tube of $\widehat{\beta}$ in
$\eta(O(\widehat{\beta},\rho))$ which is filled with the nontrivial braid
if and only if the tube of $O(\widehat{\beta},\rho)$ is. Similarly, for
$i=1,\ldots,s-1$ we may assume that $\eta^i$ sends each $\widehat{\beta}$-tube
in $O(\widehat{\beta},\rho)$ to a $\widehat{\beta}$-tube in
$\eta^i(O(\widehat{\beta},\rho))$ which has the same filling in $\beta$.

Now we can use the same property as a construction recipe for $h(\rho)$:
a tube of $\rho$ in $\eta^i(O(\widehat{\beta},\rho))$ (where $i=1,\ldots,s$)
is filled in the same way as its preimage under $\eta^i$. With this
definition, $h(\rho)$ commutes with $\beta$. Finally we are ready to
define $h(\eta)$: we take the braid $\psi(\eta)$, but modify the braids
in the tubes that terminate at positions corresponding to
$O(\widehat{\beta},\rho)$ by multiplying them on the right by $\iota_*$. In
other words, the braid $h(\eta)$ is obtained from $\eta$ as follows: we fill
those tubes of $\eta$ which connect points in $\eta^i(O(\widehat{\beta},\rho))$
to points in $\eta^{i+1}(O(\widehat{\beta},\rho))$ (with $i=0,\ldots,s-2$)
with the trivial braid, and we fill the tubes that start in
$\eta^{s-1}(O(\widehat{\beta},\rho))$ and terminate in
$O(\widehat{\beta},\rho)$ with the interior braids of $\iota_*$.
By construction, this braid $h(\eta)$ commutes with both $\beta$ and $h(\rho)$.
This concludes the proof of proposition \ref{P:section_pA}.
\end{proof}

\begin{cor}
Suppose that $\beta$ is a non-periodic reducible braid in regular form. Then
the exact sequence
$$
1\longrightarrow Z(\beta_{[1]})\times \cdots \times Z(\beta_{[t]})
\stackrel{g}{\longrightarrow} Z(\beta)  \stackrel{p}{\longrightarrow}
Z_0(\widehat{\beta}) \longrightarrow 1
$$
splits. That is,
$Z(\beta)\cong\left(Z(\beta_{[1]})\times \cdots
\times Z(\beta_{[t]})\right) \rtimes Z_0(\widehat{\beta})$.
\end{cor}

\begin{proof}
 Since $\widehat{\beta}$ cannot be reducible, the result is a direct
consequence of propositions \ref{P:section_trivial}, \ref{P:section_periodic}
and \ref{P:section_pA}.
\end{proof}


\subsection{Structure of $Z_0(\widehat{\beta})$}\label{S:C_0beta}

The proof of theorem~\ref{T:main1} is now completed by the following result.
\begin{prop}\label{P:C_0beta}
Suppose that $\beta$ is a non-periodic reducible braid, and that its tubular
braid $\widehat{\beta}$ has $m$ strands. Then $Z_0(\widehat{\beta})$ is
isomorphic either to $\Z^2$ or to a mixed braid group on $k$ strands,
where $k\leqslant m$.
\end{prop}

\begin{proof}
As usual, there are three subcases, depending wether $\widehat{\beta}$ is
trivial, periodic or \pA.
Recall that we are assuming that $\beta$ is in regular form.

Suppose first that $\widehat{\beta}=1$. In this case, $Z(\widehat{\beta})=B_m$.
Hence $Z_0(\widehat{\beta})$ contains any braid whose permutation is
consistent with $\beta$. Denote by ${\cal P}$ the following partition of
$\{P_1,\ldots,P_m\}=\{C_{1,1},\ldots,C_{m,1}\}$: we say that $P_i$ and $P_j$
belong to the same coset of ${\cal P}$ if and only if
$\beta_{[i]}=\beta_{[j]}$. By definition, a braid's permutation is
consistent with $\beta$ if and only if it preserves ${\cal P}$. Therefore,
$Z_0(\widehat{\beta})=B_{\cal P}$, and we are done. (In this case, we have
$k=m$).

If $\widehat{\beta}$ is {\pA}, it is shown in proposition~\ref{P:section_pA}
that $Z_0(\widehat{\beta})\simeq \Z^2$, so this case is already known.

Finally, suppose that $\widehat{\beta}$ is periodic. If it is a power of
$\Delta^2$, then its centralizer is the whole $B_m$, and its corresponding
permutation is trivial, so this case is equivalent to the first one.

If $\widehat{\beta}$ is periodic but not a power of $\Delta^2$, then we know
by theorems~\ref{T:centrdelta} and \ref{T:centrdelta2} that
$Z(\widehat{\beta})\simeq B_d(D_*)$, for some $d\geqslant 1$, where $D^*$
is the once punctured disk. But every base point $Q_i$ in $D^*$ corresponds
to an orbit ${\cal C}_i$ of $\widehat{\beta}$ (see figure~\ref{F:theta} in
section~\ref{S:per}), so we can define the following partition ${\cal P}'$
of $\{Q_1,\ldots,Q_d\}$: $Q_i$ and $Q_j$ belong to the same coset if and
only if $\beta_{[i]}=\beta_{[j]}$. This partition lifts by $\theta^{-1}$ to
a partition of $\{P_1,\ldots,P_m\}$, in such a way that any braid in
$B_d(D_*)$ preserves ${\cal P}'$ if and only if its corresponding permutation
in $Z(\widehat{\beta})$ is consistent with $\beta$. Therefore,
$Z_0(\widehat{\beta})\simeq B_{{\cal P}'}(D^*)$. Now it suffices to consider
the central puncture of $D^*$ as another base point, $Q_{d+1}$, and to
notice that $B_{{\cal P}'}(D^*)\cong B_{\cal P}$, where
${\cal P}={\cal P}'\cup \{\{Q_{d+1}\}\}$. To summarize, in this case we have
$Z_0(\widehat{\beta})\cong B_{{\cal P}'}(D^*)\cong B_{\cal P}$, and the
partition $\cal P$ has $k=d+1$ cosets. Since $d$ must be a proper divisor
of $m$, we get that $k=d+1<m$, and the result follows.
\end{proof}

In particular, $Z_0(\widehat{\beta})$ is isomorphic either to $\Z^2$ or to
a mixed braid group. Theorem~\ref{T:main1} is thus proven.


\section{An upper bound for the number of generators}\label{S:upperbound}

Once decomposed $Z(\beta)$, if $\beta$ is reducible, as a semi-direct product
of $(Z(\beta_{[1]})\times \cdots \times Z(\beta_{[t]}))$ and
$Z_0(\widehat{\beta})\subset Z(\widehat{\beta})$, we will define a small set
of generators for $Z(\beta)$. We will proceed by induction on the number of
strings, but we need to define first a generating set for
$Z_0(\widehat{\beta})$. We do it as follows:

\begin{prop}\label{P:genC_0}
 Let $\beta\in B_n$ be a non-periodic reducible braid, and let
$\widehat{\beta}\in B_m$ be its corresponding tubular braid. Then
$Z_0(\widehat{\beta})$ can be generated by at most $\frac{m(m-1)}{2}$ elements.
\end{prop}

\begin{proof}
 If $m=2$ then $Z_0(\widehat{\beta})$ is cyclic, so let us assume that
$m\geqslant 3$. We know by subsection~\ref{S:C_0beta} that
$Z_0(\widehat{\beta})$ is either isomorphic to $\Z^2$ or to a mixed braid
group. The case $\Z^2$ satisfies our result, so we will assume that
$Z_0(\widehat{\beta})$ is isomorphic to a mixed braid group on $k$ strings.

Mixed braid groups have been studied in~\cite{Manfredini}, where a
presentation in terms of generators and relations is given. Since we are
mainly interested in the generators, we will extract from those
in~\cite{Manfredini} a small generating set: Let ${\cal P}$ be a partition
of the set $\{1,\ldots,k\}$, having $d$ cosets of length $m_i$
(for $i=1,\ldots,d$). A generating set for $B_{\cal P}$ is given by the
following:
\begin{enumerate}

 \item For $i=1,\ldots, d$, a generating set for $B_{m_i}$ (if $m_i>1$).

 \item A generating set for the pure braid group $P_d$.

\end{enumerate}
It is clear that the first kind of generators corresponds to the movements
of the points inside a coset, while the second one corresponds to the
movement of the points of a coset with respect to those of the others.
For instance, if $k=6$ and ${\cal P}=\{\{1\},\{2,3\},\{4,5,6\}\}$, then one
possible generating set would be:
$$
   \{ \sigma_2\}\cup \{\sigma_4,\;\sigma_5\} \cup \{\sigma_1^2, \;
   \sigma_1\sigma_2\sigma_3^2\sigma_2^{-1} \sigma_1^{-1},\; \sigma_3^2\}.
$$
In order to minimise these generators we recall that $B_2$ is cyclic and,
if $m>2$, then $B_m$ can be generated by two elements. Hence, if we denote
$e_i=m_i-1$ if $m_i< 3$ and $e_i=2$ otherwise, then $e_i$ is a minimal number
of generators for $B_{m_i}$. On the other hand, a minimal number of generators
for $P_d$ is $\frac{d(d-1)}{2}$. Therefore, the minimal number of generators
for $B_{\cal P}$ is:
$$
   g_{\cal P}=\left(\sum_{i=1}^d{e_i}\right)+
   \frac{d(d-1)}{2}\leqslant  \left(\sum_{i=1}^d(m_i-1)\right)+\frac{d(d-1)}{2}
$$
$$
   = k-d+ \frac{d(d-1)}{2}=k+\frac{d(d-3)}{2} \leqslant  k+\frac{k(k-3)}{2}
   =\frac{k(k-1)}{2}.
$$
Notice that if ${\cal P}=\{\{1\},\{2\},\ldots,\{k\}\}$ (so $d=k$), then
$g_{\cal P}=\frac{k(k-1)}{2}$, and this is the worst possibility by the above
formula.

Finally we recall from proposition \ref{P:C_0beta} that $k\leqslant m$, so
that $g_{\cal P}\leqslant \frac{m(m-1)}{2}$.
\end{proof}

The first generating set $G'$ of $Z(\beta)$ that we will present
is the following: if $\beta$ is periodic or {\pA}, we have already defined
in sections~\ref{S:per} and \ref{S:pA} a minimal generating set of $Z(\beta)$,
having one or two elements. So suppose that $\beta$ is reducible. Then, by
induction on the number of strings, and by proposition~\ref{P:genC_0}, we
can suppose that we have defined $G_1,\ldots, G_t$ and $G_0$, generating sets
for $Z(\beta_{[1]}), \ldots, Z(\beta_{[t]})$ and $Z_0(\widehat{\beta})$
respectively (if some $\beta_{[i]}$ has one string, then $G_i=\emptyset$).
Then we define $G'=g_1(G_1)\cup \cdots \cup g_t(G_t) \cup h(G_0)$, which is
clearly a generating set for $Z(\beta)$.

\begin{proof}[Proof of theorem~\ref{T:main}.]
Denote $p(n)$ the upper bound proposed in theorem~\ref{T:main}, that is,
$p(n)=\frac{k(k+1)}{2}$ if $n=2k$ or $p(n)=\frac{k(k+3)}{2}$ if $n=2k+1$.
We will show that the generating set $G'$ defined above has at most
$p(n)$ elements. The case $n=2$ is trivial, so we can suppose that $n>2$ and
that the result is true for any smaller number of strings. We can also assume
that $\beta$ is non-periodic and reducible.

The strategy now is to successively replace $\beta$ by different braids,
in such a way that during each replacement step the number of generators
of its centralizer, as given by the above construction, increases.

The first modification of $\beta$ will be to replace the tubular braid
$\widehat{\beta}$ by the trivial braid. At the same time, we shall modify the
interior braids, with the aim of rendering them pairwise non-conjugate.
More precisely, we notice that, for any braid $\alpha$ with at least two
strings, the number of generators of $Z(\alpha)$ and $Z(\Delta^{2p}\alpha)$
is the same, while $\Delta^{2p}\alpha$ and $\Delta^{2q}\alpha$ are conjugate
if and only if $p=q$. Thus after multyplying each interior braid $\beta_{[i]}$
by a suitable power of twists $\Delta_{(m_i)}^2$, we can assume that all the
interior braids with at least two strings are pairwise non-conjugate, so that
$t=m$. As seen in the proof of proposition~\ref{P:genC_0}, this first
replacement has increased (or left unchanged) the number of generators of
$G_0$, according to our construction.


Suppose, without loss of generality, that $m_1=m_2=\cdots=m_d=1$,  that $m_i=2s_i$ for $i=d+1,\ldots,d+u$,
and that $m_i=2s_i+1$, for $i=d+u+1,\ldots,d+u+v$, where $d+u+v=m$. Hence $u$ is the number of interior braids
with an even number of strings, and $v$ is the number of interior braids with an odd (but greater than one) number
of strings. If $d\geqslant 2$, then we shall make further modifications to the braid $\beta$, with the aim
of lowering $d$. More precisely, if $d\leqslant 2$, then we can decrease $d$ by multiplying $\beta$ by
$\sigma_1^p$ for some $p$, where $p$ is chosen in such a way that no other interior braid of $\beta$ equals
$\sigma_1^p$. This replacement increases $u$ by one,
and decreases $d$ by two. Thus the number of generators in $G_0$ decreases by one (if $d=2$) or increases
(if $d>2$). But we would have a new interior braid, $\sigma_1^p$, yielding one new generator. Hence, the total
number of elements in $|G'|$ will not decrease. In other words, without decreasing the number of elements
of $|G'|$ we can replace $\beta$ by a braid with $d\leqslant 1$.


Denote $a=s_{d+1}+\cdots +s_{d+u}$, $\; b=s_{d+u+1}+\cdots+s_{m}$ and $S=a+b$.
Then one has $n=d+2S+v$. By induction on the number of strings, we have the following bound on the number of
elements in $G'$:
\begin{eqnarray*}
|G'| & \leqslant & \sum_{i=d+1}^{m}{p(m_i)}+ \frac{m(m-1)}{2}  \\
   & = & \sum_{i=d+1}^{d+u}{\frac{s_i(s_i+1)}{2}}+\sum_{i=d+u+1}^m{\frac{s_i(s_i+3)}{2}}+
\left( \begin{array}{c} m \\2 \end{array}\right) \\
  & = &  \sum_{i=d+1}^m{\left( \begin{array}{c} s_i+1 \\2 \end{array}\right)} + \sum_{i=d+u+1}^m{s_i} +
\left( \begin{array}{c} m \\2 \end{array}\right) \\
  & = &  \sum_{i=d+1}^m{\left( \begin{array}{c} s_i+1 \\2 \end{array}\right)} + b +
\left( \begin{array}{c} m \\2 \end{array}\right) \\
\end{eqnarray*}
where $s_i\geqslant 1$ for $i=d+1,\ldots,m$.

 Given two positive integers $x$ and $y$, one has:
$$
  \left( \begin{array}{c} x+1 \\ 2 \end{array}\right) + \left( \begin{array}{c} y+1 \\ 2 \end{array}\right) =
 \left( \begin{array}{c} x+y+1 \\ 2 \end{array}\right) - xy.
$$
This yields:
\begin{eqnarray*}
 |G'| & \leqslant & \left( \begin{array}{c} S+1 \\ 2 \end{array}\right)
 -\left(\sum_{d+1\leqslant i<j \leqslant m}{s_is_j}\right) + b +
\left( \begin{array}{c} m \\ 2 \end{array}\right).
\end{eqnarray*}

 Now we distinguish two cases. If $d=0$, then $m=u+v$ and $n=2S+v$. Also,
\begin{eqnarray*}
 |G'| & \leqslant & \left( \begin{array}{c} S+1 \\ 2 \end{array}\right)
 -\left(\sum_{1\leqslant i<j\leqslant m}{s_is_j}\right) + b +
\left( \begin{array}{c} m \\ 2 \end{array}\right) \\   & \leqslant & \left(
\begin{array}{c} S+1 \\ 2 \end{array}\right)
 -  \left( \begin{array}{c} m \\ 2 \end{array}\right) + b +
\left( \begin{array}{c} m \\ 2 \end{array}\right)
\quad = \quad \frac{S(S+1)}{2} + b.
\end{eqnarray*}

If $v=0$ one has $b=0$, so $S=a$ and
$|G'|\leqslant \frac{S(S+1)}{2}=\frac{a(a+1)}{2}$; but also $n=2k=2a$,
so $p(n)=\frac{a(a+1)}{2}$ and we are done.

If $v=1$ then $n=2S+1$, hence $k=S$ and $p(n)=\frac{S(S+3)}{2}$. But in this case
$|G'|\leqslant \frac{S(S+1)}{2}+b\leqslant \frac{S(S+1)}{2}+S = \frac{S(S+3)}{2}=p(n)$.

If $v\geqslant 2$, since $n=2S+v$ one has $k\geqslant S+1$. Then $|G'|\leqslant \frac{S(S+1)}{2} + b
< \frac{S(S+1)}{2} + (S+1) =\frac{(S+2)(S+1)}{2} \leqslant \frac{k(k+1)}{2} \leqslant p(n)$.

Therefore, the result is true if $d=0$. Suppose now that $d=1$. In this case $m=u+v+1$ and
$n=2S+v+1$. Then one has:
\begin{eqnarray*}
 |G'| & \leqslant & \left( \begin{array}{c} S+1 \\ 2 \end{array}\right)
 -\left(\sum_{2\leqslant i<j\leqslant m}{s_is_j}\right) + b +
\left( \begin{array}{c} m \\ 2 \end{array}\right) \\
& \leqslant & \left( \begin{array}{c} S+1 \\ 2 \end{array}\right)
 -  \left( \begin{array}{c} m-1 \\ 2 \end{array}\right) + b +
\left( \begin{array}{c} m \\ 2 \end{array}\right) \\
 & = & \frac{S(S+1)}{2} + b + m-1 \quad = \quad \frac{S(S+1)}{2} + b+u+v \\
 & \leqslant & \frac{S(S+1)}{2} + S+ v  \quad = \quad  \frac{S(S+3)}{2} + v.
\end{eqnarray*}

If $v=0$ then $b=0$ and $\;k=S$, so $|G'|\leqslant \frac{S(S+3)}{2} = p(n)$.

If $v=1$ then $n=2S+2$ and $k=S+1$. Then $|G'|\leqslant \frac{S(S+3)}{2}+1=
\frac{(S+1)(S+2)}{2}=p(n)$.

If $v=2$ then $n=2S+3$ and $k=S+1$. Then $|G'|\leqslant \frac{S(S+3)}{2}+2=
\frac{S^2+3S+4}{2} <\frac{(S+1)(S+4)}{2}=p(n)$.

Finally, if $v\geqslant 3$ then $n=2S+v+1$ so $k\geqslant S+v/2$. Hence
$$
  p(n)\; \geqslant \; \frac{(S+v/2)(S+v/2+1)}{2} \;=
  \;\frac{S^2+(v+1)S+v(v+2)/4}{2}
$$
$$
  \geqslant \; \frac{S(S+3)}{2} + S/2 + v/2
  \; > \; \frac{S(S+3)}{2}+ v \; \geqslant \; |G'|.
$$
Therefore, in every case $|G'|\leqslant p(n)$, and theorem~\ref{T:main} is
proved.
\end{proof}

Recall that, in example~\ref{E:purebraid}, 
we defined braids of any number of strands whose centralizer could not be
generated by less than $p(n)$ elements. Therefore, the bound given by
theorem~\ref{T:main} is the best possible one.


\section{Small generating sets}\label{S:minimal}

We saw in the previous section an upper bound for the number of generators of the
centralizer of a braid $\beta$, in terms of its number of strings. But one could obtain a
better bound if more information about $\beta$ is given. In this section we will define a
new generating set $G$ for $Z(\beta)$, which is in most cases smaller than the set $G'$
defined before. It is also the smallest possible ``natural'' generating set, in the sense
that each generator belongs to one of the $t+1$ factors in the semidirect product
decomposition in theorem \ref{T:main1}(c). Thus in a philosophical sense, $G$ is the
``right'' generating set, even though it is not in general the smallest possible one, as we
shall see at the end of this section.

If $\beta$ is periodic or {\pA}, we already know a minimal generating set, with at most two
elements. We also know a minimal generating set for any mixed braid group (see the proof of
proposition~\ref{P:genC_0}). Hence we can define $G$ by induction on the number of strands,
when $\beta$ is a reducible, non-periodic braid. We can also suppose that $\beta$ is in
regular form. We recall that the interior braids are denoted $\beta_{[1]},\ldots,
\beta_{[t]}$, and the tubular braid $\widehat{\beta}$.

\begin{definition}
We will say that $i,j\in \{1,\ldots, t\}$ are {\em permutable}
if there exists some $\eta\in Z_0(\widehat{\beta})$ such that
$\eta({\cal C}_i)={\cal C}_j$.
\end{definition}

Remark that permutability is an equivalence relation,
and the definition of $Z_0(\widehat\beta )$ says that if $i$ and $j$
are permutable then $\beta_{[i]}=\beta_{[j]}$.

Let then $\{i_1,\ldots,i_r\}\subset\{1,\ldots,t\}$ be coset representatives
for permutability. Let $G_{i_k}$ be a minimal set of generators for
$Z(\beta_{[i_k]})$, and $G_0$ be a minimal set of generators for
$Z_0(\widehat{\beta})$. Then we define
$G=g_{i_1}(G_{i_1})\cup \cdots \cup g_{i_r}(G_{i_r}) \cup h(G_H)$.
Notice that $G\subset G'$, and they coincide if and only if there is no
pair of permutable indices.

\begin{prop}
$G$ is a generating set of $Z(\beta)$.
\end{prop}

\begin{proof}
From the exact sequence of theorem~\ref{T:exactseq} it follows that,
if $G_i$ is a set of generators for $Z(\beta_{[i]})$, then a set of
generators for $Z(\beta)$ is
$G'=g_{1}(G_{1})\cup \cdots \cup g_{t}(G_{t}) \cup h(G_0)$. Hence, we just
need to show that if $j\in\{1,\ldots, t\}\backslash \{i_1,\ldots,i_r\}$,
then every element in $g_j(G_j)$ can be written as a product of elements
in $G$.

Take then $j$ as above. There must be some $i_k$ permutable with $j$, so $\beta_{[j]}=\beta_{[i_k]}$
and there is some $\eta\in Z_0(\widehat{\beta})$ such that $\eta({\cal C}_{i_k})={\cal C}_{j}$.
Notice that $G_j$ is a set of generators for $Z(\beta_{[j]})=Z(\beta_{[i_k]})$, so every $\gamma\in G_j$
can be written as a product of elements in $G_{i_k}$. Hence the braid
$\alpha=h(\eta)^{-1} g_{i_k}(\gamma) h(\eta)$ can be written as a product of elements in $G$.
Moreover, one has $\widehat{\alpha}=\widehat{h(\eta)}^{-1} 1 \widehat{h(\eta)} = 1$, and the
only nontrivial interior braids in $\alpha$ are those corresponding to ${\cal C}_j$. Since the
interior braids $h(\eta)_{i_k,l}$ for every $l$  are just powers of $\beta_{[i_k]}=\beta_{[j]}$,
and $\gamma$ commutes with $\beta_{[j]}$, it follows that for every $l$, $\alpha_{j,l}=\gamma$.
Therefore $\alpha=g_j(\gamma)$, so every element in $g_j(G_j)$ can be written as a product of elements
in $G$, thus $G$ is a generating set for $Z(\beta)$.
\end{proof}

The generating set we have just defined is, unfortunately, not always the
smallest possible one:

\begin{example}\rm
Consider the five string braid
$\beta = \sigma_3 \sigma_4 \sigma_2 \sigma_3 \sigma_1
\sigma_2 \sigma_2 \sigma_3 \sigma_4 \sigma_1 \sigma_2 \sigma_3$ --
the canonical reduction system of this braid has two round circles,
one containing punctures number 1, 2 and 3, the other punctures number
4 and 5; the tubular braid is just a full twist of the two fat
strings: $\widehat{\beta} = \sigma_1^2$. Moreover, the interior braids of
each tube is trivial. According to theorem \ref{T:main1}, the centralizer
of this braid is
$$
Z(\beta)\cong (B_3 \times B_2) \rtimes PB_2 \cong (B_3 \times \Z)\rtimes \Z
$$
and the generating set constructed in this section has four elements:
two for $B_3$, and one for each factor $\Z$. We now claim that this
generating set is not as small as possible.

Indeed, $B_3 \times \Z$ can be generated by only two elements (and thus
$Z(\beta)$ can be generated by three elements). To see this, recall that
the 3-string braid group is isomorphic to the group of the $(2,3)$-torus knot.
Thus $B_3$ has a presentation $\langle y,z\ |\ y^3z^{-2} = 1\rangle$
(with $y=\sigma_1\sigma_2$ and $z=\sigma_1 \sigma_2 \sigma_1$).
Moreover, the factor $\Z$ is generated by $\sigma_4$. Now the two generators
$(y,\sigma_4)$ and $(z,\sigma_4)$ generate $B_3\times \Z$, because
$(1,\sigma_4)$ can be written as $(y,\sigma_4)^3 (z,\sigma_4)^{-2}$.
\end{example}


\section{Some algorithmic aspects}\label{S:alg}

The aim of this section is to present the essential ingredients for an
algorithm which, for any given braid, finds a generating set of its
centralizer subgroup that matches the description of the previous sections.
Since, for any braid $\beta$ and any $k\in \Z$, the centralizer subgroups of
$\beta$ and $\beta \Delta^{2k}$ coincide, we can always assume that $\beta$
is positive.

We start by mentioning that algorithms that perform the Nielsen-Thurston
classification, and give the invariant folitations in the {\pA} case (in the
form of train tracks), are available -- notably, there are Bestvina-Haendel's
\cite{BestvinaHaendel} and of Los' \cite{Los} algorithms; and computer
implementations are available on the web \cite{Brinkmann, Hall}.

We recall briefly the idea of the two automatic structures on braid groups
that are relevant for us: for the first one, given by Garside \cite{Garside}
and Thurston \cite{ThBnaut} (and refined by El-Rifai and Morton \cite{EM}),
we think of $D_n$ has having the $n$ punctures lined up on the real line
in the disk $D$. For the second one, given by Birman, Ko, and Lee \cite{BKL},
we think of $B_n$ as having the $n$ punctures regularly spaced on the circle
of radius $1$. Apart from that, the structures are exactly analogue.
In the Garside-Thurston structure, there is a canonical way to write $\beta$
as a product of divisors of $\Delta$, namely by pushing each crossing between
two strings into a factor as far to the  left as possible. This normal form
is called the {\sl left greedy normal form}. For instance, in this normal
form all factors which are {\it equal} to $\Delta$ (not just divisors of it)
are grouped together at the very left of the product decomposition.
Analogously, Birman-Ko-Lee write
each braid as a product of divisors of $\delta$ in a left-greedy way.
If $\beta$ is a positive braid, then its {\sl super summit set} is the
subset of all elements $\alpha$ of its conjugacy class which satisfy the
following conditions:
\begin{enumerate}
\item[(i)] $\alpha$ is positive,
\item[(ii)] the writing of $\alpha$ in left greedy normal form has as
few factors as possible among all elements satisfying (i),
\item[(iii)]  the writing of of $\alpha$ in left greedy normal form has
as many factors on the left as possible equal to $\Delta$ (or $\delta$),
among all elements satisfying (i) and (ii).
\end{enumerate}
Two positive elements of $B_n$ are conjugate if and only if their super
summit sets coincide. Given $\beta\in B_n$ there is an algorithm, given
in~\cite{FrGM_conj} (which is an improvement of the algorithm in~\cite{EM}),
to compute its super summit set. It is as follows: first we repeatedly
{\sl cycle} $\beta$ (i.e.~move the first factor different from $\Delta$,
respectively $\delta$, to the end and calculate the left greedy form of the
resulting braid), until this process runs into a loop. At this point we are
guaranteed to have achieved condition (ii) above. Then we repeatedly
{\sl decycle} (i.e.~move the last factor to the front
and calculate the left greedy form of the resulting braid) until we run
into a loop. Then all elements of this loop belong to the super summit set.
Afterwards, all other elements of the super summit set can be found
recursively by conjugating already known elements by (suitable) divisors of
$\Delta$ (respectively $\delta$), and retaining the result if it belongs to
the super summit set.

This algorithm for computing the super summit set is necessary for our
purposes. Now suppose we are given a braid $\beta\in B_n$ and we want to
compute its centralizer. First we need to determine if $\beta$ is periodic,
reducible or {\pA}, and then we can use the results in this paper.

\begin{remark}
Very recently, V. Gebhardt~\cite{Gebhardt} presented a better algorithm for
the conjugacy problem in braid groups. He defined the {\em ultra summit set},
which is in general much smaller than the super summit set described here.
\end{remark}


\sh{Periodic elements}

Deciding whether a given element $\beta$ of $B_n$ is periodic is very
easy: one calculates the $n-1$st and the $n$th power of $\beta$. Then $\beta$
is periodic if and only if one of the two results is a power of $\Delta^2$.

If $\beta^{n-1}=\Delta^{2k}$ for some $k\in \N$, then $\beta$ is conjugate
to $\gamma_{(n)}^k$ (as can be easily seen from lemma~\ref{L:per_gamma_delta}),
and a conjugating element can be found explicitly using either of the
two standard algorithms. Similarly, if $\beta^n=\Delta^{2k}$, then $\beta$
is conjugate to $\delta_{(n)}^k$, and either algorithm yields an
explicit conjugating element. In either case, one can find explicitly
a generating set of the centralizer subgroup with only two elements, using
propositions \ref{P:centrdelta'} or \ref{P:centrgamma}.


\sh{Finding reducing curves of reducible elements}

After establishing that an element $\beta$ of $B_n$ is not periodic, we need
to check whether it is reducible, and if it is, we want to find explicitly
an invariant multicurve. This is, in fact, a standard part of
Bestvina-Haendel's \cite{BestvinaHaendel} and of Los' \cite{Los}
algorithms.

We want to point out one particulary elegant alternative, which is due to
Benardete, Gutierrez and Nitecki \cite{BGN2} (see also~\cite{BGN}). We think of
$D_n$ as having the $n$ punctures lined up horizontally, and we look at
Garside-Thurston's left greedy normal form. The key observation from
\cite{BGN2} is the following: suppose that $C$ is an invariant multicurve of
a braid $\beta$, and that the normal form of $\beta$ is
$\beta=\beta_1\cdot\ldots\cdot \beta_k$, where $\beta_1,\ldots,\beta_k\in B_n$
are divisors
of $\Delta$. Moreover, suppose that all components of $C$ are {\sl round}
(i.e.~actual geometric circles in $D_n$). Then we have not only that
$\beta_1\cdot\ldots\cdot\beta_n(C)=C$, but also that all components of all the
multicurves $\beta_1\cdot\ldots\cdot\beta_i(C)$ are round for $i=1,\ldots,k$.

As remarked in \cite{BGN2} this implies as a corollary that invariant
multicurves are visible as round curves in the super summit set of $\beta$,
and in particular the reducibilty of a braid is easily detectable from the
super summit set. To prove the corollary we note that $\beta$ has a conjugate
in which all components of the curve system $C$ are round; moreover,
$\beta$ and its conjugate have the same super summit set. Now cycling and
decycling this conjugate does not change the fact that there  is a {\it round}
invariant curve system, by the key observation above. At the end of the
cycling/decycling procedure we have found elements of the super summit set
which contain the desired round invariant curves.

Now it is shown in~\cite{BGN2} how to determine if a
given braid preserves a system of disjoint round curves. And there is a finite
number of these systems. Moreover, since for each element of the super summit
set we know how it can be conjugated to obtain $\beta$, we can find explicitly
all curves that belong to a reduction system for $\beta$. We can then easily
determine, by its definition, which of these curves belong to the canonical
reduction system of $\beta$. That is, we can compute the canonical reduction
system of $\beta$.

By the results in this paper, $Z(\beta)$ is then a semi-direct product of two
groups that can be computed by induction on the number of strings. Hence, it
only remains to study the case when $\beta$ is {\pA}.


\sh{Pseudo-Anosov elements: commutation with $\delta_{(n)}^k$}

Suppose that our braid $\beta$ fails the tests of periodicity and reducibility,
hence it is known to be {\pA}. We need to check if it commutes with a
periodic braids other than powers of $\Delta^2$.

We shall think of $D_n$ as having its $n$ punctures uniformly distributed
over the circle of radius $1$, and we consider Birman-Ko-Lee's
left-greedy normal form. We want to decide
algorithmically whether $\beta$ is conjugate to a braid $\alpha$ with
the property that $\alpha$ commutes with $\delta_{(n)}^k$ for some
positive integer $k<n$. If it is, we want to know the conjugating braid
explicitly. The following result yields such an algorithm.
\begin{prop}\label{P:commdelta}
Suppose that a {\pA} braid $\beta$ has a conjugate which commutes with
$\delta_{(n)}^k$ for some integer $k$. Then there exists an element $\alpha$
of the super summit set of $\beta$ which has the property that $\alpha$,
and in fact every factor of the left greedy normal form of $\alpha$,
commutes with $\delta_{(n)}^k$.
\end{prop}
\begin{proof} Let $\beta'$ be a conjugate of $\beta$ which commutes with
$\delta_{(n)}^k$. If $\beta'=\beta'_1\cdot\ldots\cdot\beta'_r$ is the
left-greedy normal form of $\beta'$, then each factor $\beta'_i$ is a divisor
of $\delta_{(n)}$ which is $\frac{2\pi k}{n}$-symmetric. This follows from
the fact that the very definition of the left-greedy normal form is
completely rotation symmetric. More precisely, the fact that two
consecutive factors $\beta'_i \beta'_{i+1}$ determine a left-greedy normal
form is not modified by rotating them. Hence, the product
$(\delta_{(n)}^{-k} \beta'_1 \delta_{(n)}^k ) \cdots
(\delta_{(n)}^{-k} \beta'_r\delta_{(n)}^k)$ is in left-greedy normal form.
Since this product equals $\delta_{(n)}^{-k} \beta'\delta_{(n)}^k=\beta'$,
whose left-greedy normal form is $\beta'_1\cdots \beta'_k$, we obtain that
$\delta_{(n)}^{-1} \beta'_i \delta_{(n)}=\beta'_i$, for $i=1,\ldots,r$.

Using the same argument inductively, we see that the cycling and decycling
procedure only ever creates braids in left greedy normal form in which all
factors are $\frac{2\pi k}{n}$-symmetric.
\end{proof}

Now we notice that it is very easy to decide if a given divisor of $\delta$ (in
the Birman-Ko-Lee context) is invariant under a given rotation. Hence one
can determine if a braid commutes with an (explicitly computable) conjugate of
$\delta_{(n)}^k$ by looking at the elements of its super summit set.


\sh{Pseudo-Anosov elements: commutation with $\gamma_{(n)}^k$}

 Now we want to determine if a given {\pA} braid commutes with a conjugate
of $\gamma_{(n)}^k$, for a given positive integer $k<n-1$. This is only
possible if there is some index $i\in \{1,\ldots,n\}$ such that $\beta$
preserves $P_i$, as can be easily seen by looking at the corresponding
permutations.

Call ${\cal P}_i=\{\{P_i\},\{P_1,\ldots,P_{i-1},P_{i+1},\ldots,P_n\}\}$, a
partition of $\{P_1,\ldots,P_n\}$. Then $\beta$ should belong to
$B_{{\cal P}_i}$. There is a natural map
$f_i\co B_{{\cal P}_i} \rightarrow B_{n-1}$ which consists of forgetting
the $i$th string. Notice that, if a braid $\alpha$ commutes with
$\gamma_{(n)}^k$ (where $P_1$ is considered to be the central point of
$D_{(n)}$) then $f_1(\alpha)$ commutes with
$f_1(\gamma_{(n)}^k)=\delta_{(n-1)}^k$.

Hence we have a necessary condition that must be satisfied. If $\beta$
preserves a puncture $P_i$, then we conjugate it to some $\alpha$ that
preserves $P_1$, and we test whether a conjugate of  $f_1(\alpha)$ commutes
with $\delta_{(n-1)}^k$ for some $k<n-1$. If this does not happen, for
$i=1,\ldots, n$, then no conjugate of $\beta$ commutes with $\gamma_{(n)}^k$.

This necessary condition is of course not sufficient. A sufficient
and testable condition is now given by the following result. Recall that, by
corollary~\ref{C:addstring}, there is an isomorphism
$\chi=(\bar{\theta}^*)^{-1}\theta^*$ from
$Z(\delta_{(n-1)}^k)$ to $Z(\gamma_{(n)}^k)$, given by adding a trivial
string at the centre of $D_{n-1}$. Notice that, if
$\zeta\in Z(\gamma_{(n)}^k)$, then $\chi(f_1(\zeta))=\zeta$. Then one has:

\begin{prop}
Suppose that $\alpha\in B_n$ preserves $P_1$, and $\tld\alpha=f_1(\alpha)$
commutes with $\delta_{(n-1)}^k$. Then the following two statements are
equivalent.
\begin{enumerate}
\item[(i)] $\alpha$ is conjugate to an element $\zeta$ of $B_n$
which commutes with $\gamma_{(n)}^k$, and the conjugating homeomorphism
preserves $P_1$.
\item[(ii)]  $\alpha$ is conjugate to $\chi(\tld\alpha)$.
\end{enumerate}
\end{prop}

\begin{proof}
The implication (ii)$\Rightarrow$(i) is immediate, by choosing
$\zeta := \chi(\tld\alpha)$.

For the implication (i)$\Rightarrow$(ii), we suppose that (i) holds, that is,
there is an element $\eta\in B_{{\cal P}_1}$ such that
$\eta^{-1} \alpha \eta= \zeta$, where $\zeta\in Z(\gamma_{(n)}^k)$. We can
apply $f_1$ to all these elements, denoting $\tld\eta=f_1(\eta)$ and
$\tld\zeta=f_1(\zeta)$.
This yields  $(\tld\eta)^{-1} \tld\alpha \:\tld\eta = \tld\zeta$,
where $\tld\alpha, \tld\zeta \in Z(\delta_{(n-1)}^k)$.

If we show that $\tld\eta\in Z(\delta_{(n-1)}^k)$, then we can apply $\chi$ to
all factors, obtaining
$\chi(\tld\eta)^{-1} \chi(\tld\alpha) \:\chi(\tld\eta) =
\chi(\tld\zeta)=\zeta$, hence $\chi(\tld\alpha)$ is conjugate to $\zeta$
which is conjugate to $\alpha$, and the result follows.

Let us then show that $\tld\eta$ commutes with $\delta_{(n-1)}^k$.
Notice that $\zeta$ is a {\pA} braid that commutes with $\gamma_{(n)}^k$.
Hence it preserves a projective foliation ${\cal F}_{\zeta}$, which is
invariant under a rotation by an angle of $\frac{2\pi k}{n-1}$. But in this
case $\tld\zeta$ also preserves ${\cal F}_{\zeta}$, with the same
stretch factor, hence it is also {\pA}. Since $\tld\alpha$ is conjugated
to $\tld\zeta$, then it is {\pA} as well, and we call ${\cal F}_{\tld\alpha}$
its corresponding projective foliation (which is also invariant under the
same rotation, since $\tld\alpha$ commutes with $\delta_{(n-1)}^k$).
Since $(\tld\eta)^{-1} \tld\alpha \:\tld\eta = \tld\zeta$, we have that
$\tld\eta$ sends ${\cal F}_{\tld\alpha}$ to ${\cal F}_{\zeta}$.

Now consider the braid $d=\tld\eta^{-1} \delta_{(n-1)}^k \tld\eta$.
It is conjugate to $\delta_{(n-1)}^k$, and hence periodic.
Moreover, it preserves ${\cal F}_{\zeta}$, so
it commutes with $\tld\zeta$. But the periodic elements in the centralizer
of $\tld\zeta$ form a cyclic group containing $\delta_{(n-1)}^k$, and
$\delta_{(n-1)}^k$ is the only element having exponent sum $(n-2)k$. Since
$d$ has exactly the same exponent sum, it follows that $d=\delta_{(n-1)}^k$.
Hence $\tld\eta$ commutes with $\delta_{(n-1)}$, and the result follows.
\end{proof}

An algorithm for testing whether a braid $\beta$ is conjugate to a braid
which commutes with $\gamma_{(n)}^k$ is now easy to construct: for each of
the $n$ punctures test whether the puncture is fixed by $\beta$, and whether
forgetting this puncture yields a braid which is conjugate to a braid
$\tld{\alpha}$ that commutes with $\delta_{(n)}^k$. (We know how to do
this, by the results of the previous subsection). For each puncture
that does satisfy this property, test whether $\chi(\tld\alpha)$  (which is
obtained from $\tld{\alpha}$ by adding a ``trivial'' string in the
centre), is conjugate to $\beta$. If, for one of the punctures, this is the
case, then the answer is ``yes'', otherwise ``no''.


\sh{Pseudo-Anosov elements: finding roots}

It remains to describe a last step for computing a generating set for
$Z(\beta)$, when $\beta$ is {\pA}. We assume that we have already computed
the subgroup $\langle \rho \rangle$ of periodic braids commuting with $\beta$.
Then we can multiply $\beta$ by a suitable power of $\rho$, to obtain a braid
$b$ that preserves the singular leaves of the projective foliations
corresponding to $\beta$. Then we know that
$Z(\beta)=\langle \alpha \rangle \times \langle \rho \rangle$, where $\alpha$
is the smallest possible root of $b$.

The last problem, therefore, is to determine whether a given {\pA} braid $b$
has a $k$th root, for given $k$, and to compute that root. This problem has
been solved in~\cite{Sty} (generalised to all Garside groups in \cite{Sibert}).
Moreover, since the number of possible values of $k$ is finite (we are
assuming that $b$ is positive), we have an algorithm for computing $\alpha$,
thus a generating set for $Z(\beta)$.

{\bf Acknowledgements } We are grateful to a number of people for discussions
and valuable ideas that greatly contributed to this research. The examples of
Nikolai V.\ Ivanov \cite{IvanovTalk,IvanovRecent}, which we learned about
through discussions with Mustafa Korkmaz, were an important inspiration and
greatly helped us clarify our ideas. It was thus from Ivanov (via Korkmaz)
that we learned that the number of generators may have to grow quadratically
with the number of strings, contradicting a conjecture in \cite{GMFr}. We are
very grateful to Sang Jin Lee who later, but independently of Ivanov, came up
with his examples, conjectured that they represent the worst case, and kindly
communicated these ideas to us by email. (Hessam Hamidi-Teherani found the
same examples as Lee immediately after listening to Ivanov's talk, but we
didn't learn this until very recently.) We also thank David Bessis for useful
discussions, and Joan Birman for telling us about the references~\cite{BGN}
and \cite{BGN2}.


\end{document}